\documentclass[12pt]{article} 
\usepackage{verbatim,amsmath,amssymb}
\usepackage{epsfig}
\usepackage{graphicx}
\usepackage{index}
\usepackage{amscd}
\input math111.def
\numberwithin{equation}{section}
\newtheorem{theor}{Theorem}[section]
\newtheorem{defi}[theor]{Definition}

\newtheorem{prop}[theor]{Proposition}

\def\bs{\vskip 0.4cm}
\def\supp{\mbox{\rm supp }}

\newif\ifmarglab
\makeatletter
\def
\label#1{\@bsphack\ifmarglab\marginpar{LAB:#1}\fi\if@filesw {\let\thepage\relax
   \def\protect{\noexpand\noexpand\noexpand}%
   \edef\@tempa{\write\@auxout{\string
      \newlabel{#1}{{\@currentlabel}{\thepage}}}}%
   \expandafter}\@tempa
   \if@nobreak \ifvmode\nobreak\fi\fi\fi\@esphack}
\makeatother

\begin{document}
\title{Compensated Compactness,
Separately convex Functions \\ 
and Interpolatory Estimates between Riesz Transforms and Haar Projections}
\author{ Jihoon Lee\thanks{J.L. supported by Korean Research Foundation Gr. Nr. KRF-2006-311-C00007.}, Paul F. X . M\"uller \thanks{P.F.X.M.  is supported
by the Austrian Science foundation (FWF) Pr.Nr. P15907-N08.} and Stefan M\"uller}
\date{Dezember  6 $^{{\rm th}},$ 2007.}
\maketitle

\tableofcontents
\newpage

\section{The main Results}
In  this work we prove sharp interpolatory estimates
that  exhibit a  new link between Riesz transforms and directional projections of the Haar system in $\bR^n . $
To a given direction $ \varepsilon \in \{ 0, 1 \}^n , \varepsilon \ne ( 0,\dots , 0 ) , $
we let $P^{(\varepsilon)}$ be the orthogonal projection onto the span of those Haar functions 
that oscillate along the coordinates $\{ i : \varepsilon_i = 1\} . $ 
When  $ \varepsilon_{i_0} = 1 $
the identity operator and the  Riesz transform $ R_{i_0}  $
provide a   logarithmically convex estimate for 
the  $L^p$ norm  of
$ P^{(\varepsilon)},$
see Theorem~\ref{th1}. 
Apart from its intrinsic interest  Theorem~\ref{th1} 
has direct applications to variational integrals, the theory of
compensated compactness, Young measures, and to the relation between 
rank one  and quasi convex functions. 
In particular we exploit our Theorem~\ref{th1} 
in the course of proving a conjecture of L. Tartar 
on  semi-continuity  of separately convex integrands; 
see Theorem~\ref{tatarconj}.
\subsection{Interpolatory Estimates}
We first recall the definitions of the Haar system in $\bR^n,$
indexed and supported on dyadic cubes, its associated directional Haar
projections  and the usual Riesz transforms; 
thereafter we  state the main theorem of this paper.

Let $\cD $ denote the collection of dyadic intervals in 
the real line. Thus $I \in \cD $ if there exists $i \in \bZ $ and $k \in \bZ$
so that $ I = [ i 2^k, ( i + 1 ) 2^k [ . $
Define the Haar function over the unit interval as 
$$h_{[0,1[} = 1_{[0,1/2[}-  1_{[1/2,1[}. $$
The $L^\infty$ normalized Haar system $\{h_I : I \in \cD\}$
is obtained from $h_{[0,1[}$ by rescaling. 
Let $I \in \cD,$ let $l_I $ denote the left endpoint of $I,$ 
thus   $l_I = \inf I .$ Then put
$$ h_I (x)=  h_{[0,1[}\left( \frac{ x - l_I}{|I|}\right), 
\quad\quad x \in \bR. $$
Thus defined, the 
Haar system $\{h_I : I \in \cD\}$
is a complete orthogonal system in $L^2(\bR).$ 
Next we recall its $n $ dimensional analog. 
Let $I_1 , \dots, I_n $ be dyadic intervals 
so that $ |I_i | = |I_j| ,$ where $1 \le i,j \le n . $ 
Define the dyadic cube $Q \sb \bR^n , $
$$ Q = I_1 \times \cdots \times I_n . $$
Let $\cS $ denote the collection of all dyadic cubes in  $\bR^n . $
To define the associated Haar system consider first
$ \cA = \{ \varepsilon \in \{ 0, 1 \}^n : \varepsilon \ne ( 0,\dots ,  0 ) \} . $
For $Q = I_1 \times \cdots \times I_n\in \cS $
and $  \varepsilon = (\varepsilon_1, \dots , \varepsilon_n) \in \cA $
let 
\begin{equation}\label{4maerz2} 
h_Q^{(\varepsilon)}(x) = \prod_{i = 1 }^{n}  
h_{I_i}^{\varepsilon_i} (x_i) , \quad\quad x = ( x_1 ,\dots, x_n ) . 
\end{equation}
We call $\{h_Q^{(\varepsilon)} : Q\in \cS, \varepsilon\in \cA\} $
the Haar system in $\bR^n . $
It is a complete orthogonal system in $L^2(\bR^n).$
Hence for $ u \in  L^2(\bR^n),$
\begin{equation}\label{10juli1} 
u = \sum_{\varepsilon\in \cA,\, Q\in \cS} 
\la u ,  h_Q^{(\varepsilon)}\ra  h_Q^{(\varepsilon)} |Q|^{-1}, 
\end{equation}
where the series on the right hand side converges unconditionally in 
$ L^2(\bR^n).$
For $\varepsilon\in \cA$ define the associated directional projection on $ L^2(\bR^n)$
by
$$P^{(\varepsilon)} (u) = \sum_{ Q\in \cS} 
\la u ,  h_Q^{(\varepsilon)}\ra  h_Q^{(\varepsilon)} |Q|^{-1}, \quad\quad u \in  L^2(\bR^n). $$
The operators $P^{(\varepsilon)}, \varepsilon \in \cA ,$ project 
 onto orthogonal
subspaces of $ L^2(\bR^n)$ so that 
\begin{equation}\label{14aug071}
 u = \sum_{\varepsilon\in \cA}P^{(\varepsilon)} (u) \quad{\rm and}\quad 
\|u\|_2^2 = \sum_{\varepsilon\in \cA}\|P^{(\varepsilon)} (u)\|_2^2
.
\end{equation}
Let $\cF$ denote the Fourier transformation on $\bR^n  $  given  as 
$$ \cF (u)(\xi) = \int_{\bR^n} e^{- i\la x,\xi\ra} u(x) dx, \quad \xi \in \bR^n, 
\quad x \in \bR^n. $$
The  Riesz transform 
$R_i$  $(1 \le i \le n )$ 
is a Fourier multiplier defined by    
$$
R_i(u)(x) = -\sqrt{-1} \cF^{-1}\left(\frac{\xi_i}{|\xi|}\cF(u)(\xi)\right)(x)
\quad\text{where}\quad \xi = ( \xi_1,\dots , \xi_n ).
$$
The analytic backbone of this paper is the following theorem
showing that  the norm in $ L^p(\bR^n)$
of  $P^{(\varepsilon)}(u)$ is 
dominated through  a logarithmically convex estimate  
by  $R_{i_0}(u),$ 
provided that 
 a carefully analyzed relation holds between
${i_0}$ (appearing in the Riesz transform) and  
$\varepsilon$  defining the directional projections 
$P^{(\varepsilon)}.$
\begin{theor} \label{th1}
Let $1 < p < \infty $ and 
$1/p + 1/q = 1 .$
For $ 1 \le i_0 \le n $ 
define 
$$ 
\cA_{i_0} = \{ \varepsilon \in \cA : 
 \varepsilon = (\varepsilon_1, \dots  \varepsilon_n) 
\quad\text{and}\quad \varepsilon_{i_0} = 1 
\}. $$
Let $u \in L^p(\bR^n).$
If  
 $\varepsilon \in \cA_{i_0}$
then $P^{(\varepsilon)}$ and $R_{i_0}$ are related by 
 interpolatory estimates in $L^p(\bR^n),$
$$
 ||P^{(\varepsilon)}(u)||_p \le
C_p \|u\|_p^{1/2}\|R_{i_0}(u)\|_p^{1/2}  \quad\quad\text{if} \quad  p \ge 2 ,
$$
and
$$
 ||P^{(\varepsilon)}(u)||_p \le
C_p \|u\|_p^{1/p}\|R_{i_0}(u)\|_p^{1/q} \quad \quad\text{if} \quad  p \le 2 .
$$
\end{theor}
The exponents $(1/2, 1/2)$ for $ p\ge 2$ and $
( 1/p , 1/q)$ 
for  $ p\le 2$ appearing in  Theorem~\ref{th1} are sharp. 
We show in Section~\ref{sharp} 
that for $\eta>0,$
$ 1 <p<\infty $ and $N >>1 $ there exists
$u = u_{\eta, p, N } \in L^p $ so that 
$$ ||P^{(\varepsilon)}(u)||_p \ge
N \|u\|_p^{1/2-\eta}\|R_{i_0}(u)\|_p^{1/2+\eta}  \quad\quad\text{if} 
\quad  p \ge 2 ,
$$
and 
$$  ||P^{(\varepsilon)}(u)||_p \ge
N \|u\|_p^{1/p-\eta}\|R_{i_0}(u)\|_p^{1/q+\eta} \quad \quad\text{if} 
\quad  p \le 2 .$$

\paragraph{A first  consequence of Theorem ~\ref{th1}.}
In the next subsection we will show how 
Theorem ~\ref{th1} is used  in 
problems originating in the theory of compensated compactness.
To this end we formulate here a concise inequality that follows from 
the  above interpolatory estimates, and record  
its immediate consequences.
See \eqref{27juli07.1}--\eqref{9juli3}.
 
Let $ 1 \le j \le n . $
Let $e_j \in \cA $ denote the unit vector in $\bR^n $ 
pointing in the positive direction of the $j- th$  coordinate axis,
$e_j = ( 0, \dots,1,\dots , 0) , $ where 
$1$ appears in the $j- th $ entry.
By \eqref{14aug071}
$$ u -  P^{(e_j)}  ( u ) = \sum_{\varepsilon \in \cA \sm \{ e_j \}}
P^{(\varepsilon)} ( u ) . $$
The above identity and 
the estimates of  
Theorem ~\ref{th1} combined yield the inequality
\begin{equation}
\label{9juli1}
 || u -  P^{(e_j)}  ( u )||_p \le 
C_{p,n} \|u\|_p^{1/2}\left[\sum_{\substack{1\le i \le n  \\i\ne j }} \|R_{i}(u)\|_p
\right]^{1/2}, \quad\quad p\ge 2.
\end{equation}
On $L^p ( \bR^n,  \bR^n )$ 
 define the vector valued projection $P$ by 
$$ P(v) = \left(  P^{(e_1)}  ( v_1 ), \dots ,  P^{(e_n)}  ( v_n ) \right) , $$
where $ v :  \bR^n \to   \bR^n ,$ $ v = ( v_1, \dots, v_n ). $ Applying 
\eqref{9juli1} to  each component of $v$ yields 
\begin{equation}
\label{27juli07.1}
\| v - P(v) \|_p \le C_{p,n}\|v\|^{1/2} \cdot \left( \sum_{i=1}^n \sum_{j = 1, j \ne i } ^n \|R_i(v_j)\|_p\right)^{1/2}
\end{equation}
Assume now that  
 $  ( v_{r,1}, \dots,  v_{r,n} ) $
is a sequence in 
$L^p ( \bR^n,  \bR^n)$ so that 
\begin{equation}
\label{9juli2} 
 \lim_{r \to \infty} \|R_i (  v_{r,j} )\|_p = 0 \quad\text{for}\quad 1 \le i \le n , i \ne j .
\end{equation}
The assumption \eqref{9juli2} and the estimate \eqref{27juli07.1}  imply that 
\begin{equation}
\label{9juli3}
  \lim_{r \to \infty}\|   ( v_{r,1}, \dots,  v_{r,n} ) - P\left( ( v_{r,1}, \dots,  v_{r,n} )\right) \|_p = 0. 
\end{equation}
Being able to draw the conclusion \eqref{9juli3}  from  the hypothesis 
\eqref{9juli2} 
provided the  main impetus for proving Theorem ~\ref{th1}.
\subsection{Lower semi-continuity and compensated compactness}
Here we
provide a frame of reference for the problems considered in this paper.
We review briefly some of the ideas of the  theory of 
compensated compactness which has been developed by 
F. Murat and L. Tartar \cite{mur1,mur3,tat1,tat2}. 
\paragraph{Weak lower-semicontinuity and differential constraints.}
Fix a system of first-order, linear differential operators $\cA . $
It is given by matrices $A^{(i)} \in \bR^{p\times d}, \,i \le n,$ so 
that 
$$ \cA (v) = \sum_{i = 1}^n A^{(i)} \pa_i (v) , $$
where $v : \bR ^n \to \bR^d $ and $ \pa_i$ denotes the partial differentiation
with respect to the $i-$th coordinate.  To  $ \cA $ we associate the cone 
$\Lambda \sbe \bR^d $ of  ``dangerous'' amplitudes. It consists of those  $a \in  \bR^d $
for which there is a vector of frequencies  $\xi \in  \bR ^n, \xi \ne 0 ,$
so that for any smooth $h : \bR \to \bR $ the function
$$ w(x) = a h(\la \xi,x\ra),$$
satisfies $$\cA(w)=0 . $$
Thus,  to 
 $a \in \Lambda$ there exists  a non-zero $\xi \in  \bR ^n,$
so that $\cA(w_m)=0  $ for the increasingly oscillatory sequence 
$$ w_m(x) = a \sin (m\la \xi,x\ra), \quad m \in \bN .$$ 
Since $\xi\ne0 $  there is $i_0 \le n $ so that the sequence of partial derivatives 
$ \pa _{i_0}  w_m$ is unbounded while  $\cA(w_m)=0  .$ 
In other words, the linear differential constraint  $\cA(w)=0  $ does not imply any control 
on the partial derivative $ \pa _{i_0} .$ 
Expressed formally, the cone of ``dangerous'' amplitudes is given as 
$$ \Lambda = \left\{ a \in  \bR^d : \exists \xi \in  \bR ^n \sm  \{ 0 \} \quad{ \rm such \,\,that }\quad 
\sum_{i=1}^n  \xi_i  A^{(i)} (a) = 0 \right \}.
$$

The methods of compensated compactness allow one to exploit a given set of
 information
on the differential constraints  $\cA (v)$ (respectively  on  $\Lambda$)
to analyze the limiting behaviour of non-linear integrands acting on 
 $v $ under weak conmvergence
Consider a sequence of functions $v_r : \bR^n \to \bR^d $ so that 
\begin{equation}
\label{vor1} v_r \rightharpoonup v \quad{\rm weakly \,\,in }\quad L^p(\bR^n,\bR^d) ,
\end{equation}
and
 \begin{equation}
\label{vor2}
\cA(v_r)  \quad{\rm precompact \,\, in } \quad W^{-1,p}(\bR^n,\bR^d).
\end{equation}
The following comments are included to clarify the relation  between 
 the hypotheses
 \eqref{vor1} and \eqref{vor2}.
\begin{enumerate}
\item
Had we imposed, instead of \eqref{vor1}, that 
 $v_r \to v $ strongly in  $L^p(\bR^n,\bR^d) ,$ then 
\eqref{vor2} would hold automatically. 
\item More subtle aspects of the interplay
between \eqref{vor1} and \eqref{vor2} are depending on the structure of 
$\cA $ or $\Lambda . $ For instance, in the special case when 
$\cA (v)$ controls all partial derivatives of $v,$  
we use 
Sobolev's {\it compact} embedding theorem to see that 
 \eqref{vor2}, implies  that 
 $v_r \to v $ strongly in  $L^p(\bR^n,\bR^d) .$  This case occurs when $\Lambda = \{0\}, $ 
\item
The generic (and most interesting) case arises when $\cA (v)$ fails to  control some of 
 the  partial derivatives of $v.$
This occurs when  $\Lambda \ne \{0\} .$ 
\end{enumerate}
 In the generic case one  goal of the theory is to isolate sharp  conditions on
a given $f : \bR^d \to \bR $ that {\it compensate} for the lack of 
compactness 
provided
by $\cA ,  $ and ensure that 
 \eqref{vor1} and \eqref{vor2} imply
\begin{equation}
\label{con3}
\liminf_{r\to \infty} \int_{ \bR^n} f(v_r(x))\vp(x) dx \ge \int_{ \bR^n}  f(v(x))\vp(x) dx, \quad\quad \vp 
\in C_o^+ (\bR^n).\end{equation}
Here (and below) $C_o^+ (\bR^n)$ denotes the set of non-negative compactly supported continous functions on 
$\bR^n.$
 Note that up to growth conditions on $f$ and up to 
passing to  subsequences  of $v_r ,$   the condition \eqref{con3}
states that 
$$ {\rm weak \,\,limit } \quad f(v_r)  \ge f(v) . $$
In summary, based on knowledge of $\cA$ or $\Lambda$
one goal of the theory of compensated compactness 
aims at describing and classifying those non-linearities 
$f : \bR^d \to \bR $ for which \eqref{vor1} and \eqref{vor2} imply
\eqref{con3}. 
\paragraph{Classical results on compensated compactness.}
We assume now that \eqref{vor1} and \eqref{vor2} hold
and that the differential operator $\cA$ satisfies the so called constant rank hypothesis;
for its definition see below.
The classical results of compensated compactness, as developed by 
F. Murat and L. Tartar \cite{mur1,mur3,tat1,tat2} 
assert that a general non-linearity $f$ satisfies \eqref{con3}
precisely when it is $\cA- $quasi-convex. Furthermore,
in the special case of a quadratic integrand $f(a) = \la Ma, a \ra $
the constant rank hypothesis is not needed and 
the conclusion \eqref{con3} is equivalent to $\Lambda -$convexity
of  $f(a) = \la Ma, a \ra .$ 
We state now explictely the characterizations mentioned above,
and recall the notions of $\Lambda -$convexity, 
$\cA -$quasi-convexity, and the constant rank hypothesis 
on $\cA . $

A function $ f : \bR^d \to \bR $ is  $\Lambda -$ convex if 
$$
f(\lambda a + (1-\lambda)b) \le  \lambda f(a)  +  (1-\lambda)f(b) , \quad a-b \in \Lambda,\,\, 0<\lambda < 1 . $$
The following result is due to  
F. Murat \cite{mur1}, \cite{mur2} and L. Tartar \cite{tat2}. 
\begin{prop} \label{necessary} If for every sequence 
$v_r : \bR^n \to \bR^d $ , the hypotheses   \eqref{vor1} and \eqref{vor2} 
imply \eqref{con3}, then $ f : \bR^d \to \bR $ is  $\Lambda -$convex.
\end{prop} 
Thus $\Lambda -$convexity is a necessary condition on $f$ for 
\eqref{vor1} and \eqref{vor2}  to imply \eqref{con3}.
If, moreover $f$ is quadratic, 
 $$f(a) = \la Ma , a \ra , \quad M \in \bR^{d\times d}, a \in \bR^d ,$$
then  $\Lambda -$convexity is already sufficient. This is the content of 
the following result by L. Tartar \cite{tat2}.
\begin{theor} \label{quadratic} Assume that $f$ is quadratic and $\Lambda -$convex.
Then, for every sequence 
$v_r : \bR^n \to \bR^d $ ,    \eqref{vor1} and \eqref{vor2} 
imply \eqref{con3}.
\end{theor} 

We next review the results beyond the case of quadratic integrands.
They involve the notion of  $\cA -$quasi-convexity and the constant rank hypothesis.
We define $f: \bR^d \to \bR $ to be $\cA -$quasi-convex if 
\begin{equation}
\label{aquasi}
 \int_{[0,1]^n} f(a +u(x))dx \ge f(a), 
\end{equation}
for each smooth and $[0,1]^n$ periodic  $u : \bR^n \to \bR^d ,$
that satisfies $\int_{[0,1]^n} u = 0 $ and $\cA(u) = 0.$
Note that \eqref{aquasi} asks for Jensen's inequality to hold under the decisive restriction that 
 $\cA(w) = 0.$
It was proved essentially by C.B. Morrey \cite{morr} that  $\cA -$quasi-convexity implies $\Lambda -$convexity
(see  \cite{daco1}).
The linear differential operator $\cA $ satisfies the 
constant rank hypothesis if there exists $r \le n $ so that 
$$ rk(A(\xi) ) = r ,\quad \xi \in \bS^{n-1},$$
where $$A(\xi) = \sum_{i=1}^n \xi_i A^{(i)}.$$
The next theorem 
provides a full characterization of those
integrands $f$ for which   \eqref{vor1} and \eqref{vor2} 
imply \eqref{con3}.

\begin{theor}[\cite{mur3}]\label{general}
Let $0\le f(a) \le C ( 1 + |a|^p) $ and assume that $\cA $ satisfies the constant rank hypothesis. 
Then $f: \bR^d \to \bR $ is  $\cA -$ quasi-convex if and only if \eqref{vor1} and \eqref{vor2} 
imply \eqref{con3}.
\end{theor}
A crucial component  in the  proof of Theorem~\ref{general} 
links the constant rank hypothesis and $\cA-$ quasi-convexity as follows:  
\begin{enumerate}
\item
Let $v :\bR^n \to \bR^d$ be  $[0,1]^n$ periodic
and of mean zero in $[0,1]^n.$ Under the constant rank hypothesis, 
there exists a decomposition of $v$ as
$$ v = u + w ,$$
where 
$$\cA (u) = 0 \quad{\rm  and }\quad  \|w \|_{L^p([0,1]^n)} \le C \|\cA( v)\|_{W^{-1,p}([0,1]^n)} . 
$$ 
The decomposition can be expressed in terms of an  explicit Fourier multiplier,
for which standard
$L^p $ estimates are available, provided that
the constant rank hypothesis holds.
\item
Let now  $v_r \in L^p([0,1]^n ,\bR^d)$  be a sequence of  $[0,1]^n$ periodic,
mean zero functions so that  $\cA(v_r ) \to 0 $ in $W^{-1,p}.$  Then, 
by the foregoing remark,  we may split   $v_r$ as  
$ v_r = u_r + w_r $ 
so that 
\begin{equation}\label{classic0}
\cA (u_r) = 0 \quad{\rm and }\quad \|w_r \|_p \to 0 .
\end{equation}
\item Assume moreover  that   $f$ is  $\cA-$quasi-convex. 
The decomposition 
\begin{equation} \label{classic1}
v_r = u_r + w_r 
\end{equation}
with the properties \eqref{classic0}
satisfies then 
\begin{equation} \label{classic2}
 \int_{[0,1]^n} f(a +u_r(x))dx \ge f(a), \quad{\rm and}\quad \|w_r \|_p \to 0 . 
\end{equation}
\end{enumerate}
\paragraph{Separately convex integrands.}
Wide ranging applications illustrate the power of Theorem~\ref{general}, 
yet there are 
important linear differential constraints $\cA ,$ 
for which the constant rank hypothesis does not hold and the classical proof
does not apply. Among the earliest examples
considered is the following $\cA_0, $ defined as 
$$
\left(\cA_0 (v)\right)_{i,j} =\begin{cases} \pa_i v_j & \quad i \ne j ;\\
                                          0       & \quad i = j,
                           \end{cases} 
 $$
where $v : \bR^n \to \bR^n. $ Observe that  for $ v = (v_1,\dots,v_n) $
the condition $\cA_0 (v)=0$ holds precisely when $v_i : \bR^n \to \bR $ is actually a function of the variable $x_i$ alone, that is  $v_i (x) = v_i(x_i). $ By a direct calculation, the cone of 
dangerous amplitudes associated to $\cA_0 $
is given as
$$ \Lambda_0 = \bigcup_{i=1}^n \bR e_i , $$
where $\{e_i\}$ denotes the unit vectors in $ \bR^n .$  It follows that 
the  $\Lambda_0 -$convex functions are just 
separately convex functions on $ \bR^n .$ 

For the operator $\cA_0$  the constant rank hypothesis,
does not hold, since 
$ {\rm ker} A_0 (\xi) =0 $ for $\xi \in \{e_1,\dots,e_n\}$ and 
$  {\rm ker} A_0 (e_i) =\bR e_i, \,i\le n. $
As a result the classical theory of compensated compactness
for non quadratic functionals does not apply to the operator $\cA_0 .$
Nevertheless it is an important  consequence of the interpolatory estimates in Theorem~\ref{th1}
that separately convex functions  yield  weakly semi-continuous integrands on sequences 
$v_r: \bR^n \to \bR^n $ for which  $\cA_0 (v_r)$   is  precompact in $W^{-1,p}(\bR^n,\bR^d).$
The following theorem  verifies a conjecture formulated  by  L.Tartar \cite{tat7}.

\begin{theor}\label{tatarconj}Let $1 <p <\infty. $  Assume that 
 $f: \bR^n \to \bR $  is   $\Lambda_0  -$ convex and satisfy $0\le f(a) \le C ( 1 + |a|^p) .$
Let 
$v_r : \bR^n \to \bR^n $ satisfy  
 \begin{equation}
\label{16okt0710}
v_r\rightharpoonup v \quad{\rm weakly \,\,in }\quad L^p(\bR^n,\bR^n) ,
\end{equation}
and
 \begin{equation}\label{16okt0711}
\cA_0 (v_r)  \quad{\rm precompact \,\, in } \quad W^{-1,p}(\bR^n,\bR^n).
\end{equation}
Then, 
\begin{equation}\label{5dez071}
\liminf_{r\to \infty} \int_{ \bR^n} f(v_r(x))\vp(x) dx \ge \int_{ \bR^n}  f(v(x))\vp(x) dx, \quad\quad \vp 
\in C_o^+ (\bR^n).
\end{equation}
\end{theor}
As discussed in \cite{smue2} this result 
implies that gradient Young measures supported on diagonal entries 
are laminates, and this in turn  gives an interesting relation between
rank-one convexity and quasi-convexity on subspaces with 
few rank-one directions.

In the approach of  the present paper
we fully exploit the methods introduced in \cite{smue2}. 
We base the proof of Theorem~\ref{tatarconj}
on the decomposition
given by the directional Haar projection
$$ v = P(v) + \{ v - P(v)\} ,$$
invoke the interpolatory estimates of Theorem~\ref{th1}, 
and  use the fact that $\Lambda_0-$convexity yields Jensen's inequality on the range of
$P:$ 
\begin{enumerate}
\item By  inequality \eqref{27juli07.1},  the norm of   $\{ v - P(v)\}$ 
in $L^p$ is controlled by the norm of $ \cA_0(v)$
in $W^{-1,p} . $
\item The operator $\cA_0$  does not exert any control over $P(v).$
It is  $\Lambda_0-$convexity that compensates for that.
Indeed when $f$ is separately convex we have the following form of Jensen's inequality
\begin{equation}
\label{16okt074}
 f\left(\int_{[0,1]^n} P(v) dx \right)\le \int_{[0,1]^n} f(P(v)) dx . 
\end{equation}
By rescaling of \eqref{16okt074} we get
\begin{equation}
\label{18okt071}
f( {E_M} (P(v))) \le  {E_M}( f(P(v))), \quad v \in L^p( \bR^n, \bR^n ) , \quad M \in \bZ ,
\end{equation}
where
 ${ E_M}$ denotes the conditional expectation operatpor given as 
$${ E_M} (g) (x) = \sum_{ \{R \in \cS : |R| = 2^{-Mn} \}}\int_{\bR^n} g(y) \frac{dy}{|R|} 1_R (x) , \quad g \in L^p( \bR^n).
$$ 
 
We verify  \eqref{16okt074} below. The proof is based on the observation  
that Haar functions are exactly 
localized, three-valued martingale differences.  
\item Assume  that  $ f$ is   separately convex  and that  
  $v_r \in L^p( [0,1]^n ,\bR^n)$  is a sequence
 of  $[0,1]^n$ periodic, mean zero functions so that
$\cA_0 (v_r) \to 0 $ in $W^{-1,p} . $  
With $ u_r = P(v_r)$  and  $w_r = \{ v_r - P(v_r)\} ,$ the 
decomposition 
\begin{equation}\label{ersatz1}
 v_r = u_r + w_r 
\end{equation}
 satisfies the central properties
\begin{equation} \label{ersatz2}
 \int_{[0,1]^n} f(a +u_r(x))dx \ge f(a), \quad{\rm and}\quad \|w_r \|_p \to 0.
\end{equation}
The splitting \eqref{ersatz1} with the property \eqref{ersatz2} 
is parallel to the classical decomposition \eqref{classic1} and \eqref{classic2}
 based on Fourier multipliers and the constant rank hypothesis.
\end{enumerate}
\paragraph{Jensen's inequality on the range of $P$.} We prove \eqref{16okt074}
by induction over the levels of the Haar system.
Fix $e_j ,$  the unit vector in $\bR^n $ pointing along  the $j-$th coordinate axis
and a dyadic cube  $Q= I_1\times \cdots \times I_n .$
The restriction of 
$h_Q^{(e_j)} $ to the cube $Q$ is a function of $x_j$ alone, indeed
$$  h_Q^{(e_j)} (x) = h_{I_j}(x_j), \quad x \in Q . $$
Hence for $a = (a_1,\dots, a_n)$ and $c= (c_1, \dots ,c_n) $ we have the identity
\begin{equation}\label{16okt071}
\begin{aligned}
&\int_Q f(a_1 + c_1  h_Q^{(e_1)}(x), \,\dots \,, a_n + c_n  h_Q^{(e_n)}(x) ) dx \\
  & = 
\int _Q  f(a_1 + c_1  h_{I_1}(x_1), \,\dots \,,  a_n + c_n   h_{I_n}(x_n) ) dx . 
\end{aligned}
\end{equation} 
Using \eqref{16okt071} and applying Jensen's inequality to each of the variables $x_1, \dots, x_n$  of the 
separately convex integrand $f$ gives
\begin{equation}\label{16okt073}
\begin{aligned}
  \int_Q f(a_1 + c_1  h_Q^{(e_1)}(x), \,\dots \,, a_n + c_n  h_Q^{(e_n)}(x) ) dx 
\ge |Q| f(a).
\end{aligned}
\end{equation}  
Next we fix $v = (v_1,\dots,v_n) \in L^p(\bR^n,\bR^n)$ and assume that $v_j $ is finite linear combination of 
Haar functions and not constant over the unit cube. 
Define
$$ A_{k,j} = \sum_{ \{ Q \in S : |Q| = 2^{-kn}\}} c_{Q,j}  h_Q^{(e_j)} , \quad  c_{Q,j} = \la v_j ,  h_Q^{(e_j)}\ra |Q|^{-1}.$$
Choose $M\in \bN $ and put 
$$ S_{M,j} = \sum_{k=-\infty} ^M A_{k,j} . $$
By our assumption on $v_j$ the sum defining  $S_{M,j}$ is actually finite, and there exists 
  $  M_0 $ with $M_0 \ge 0$ so that 
$$ S_{M_{0},j} = P^{(e_j)} (v_j) , \quad   1 \le j \le n . $$
Choose now  $ M \le M_0 .$ 
 Fix a dyadic cube $Q$ contained in $[0,1]^n$  with $|Q| = 2^{-Mn}. $  
Note that $  S _{M-1,j} $ is constant on $Q, $ and put $a_j  =   S _{M-1,j}(y) $ where $y \in Q  $
is chosen arbitrarily.
Furthermore, 
$$   A_{M,j}(x) = c_{Q,j}  h_Q^{(e_j)}(x) , \quad x \in Q .$$
Then, using $  S_{M,j} =  S_{M-1,j} + A_{M,j} $ and \eqref{16okt073} we obtain 
\begin{equation}\label{16okt074a}
\begin{aligned}
 \int_Q f\left( S_{M,1}(x) , \dots , S_{M,n}(x)\right) dx
&= \int _Q f\left(a_1 +  c_{Q,1}  h_Q^{(e_1)}(x) , \dots, a_n +  c_{Q,n}  h_Q^{(e_n)}(x) \right) dx \\
& \ge  |Q|  f\left( S_{M-1,1}(y) , \dots , S_{M-1,n}(y)\right) . 
\end{aligned}
\end{equation}
It follows from  \eqref{16okt074a} by  taking the sum over $ Q  \sb [0,1]^n$  with $|Q| = 2^{-Mn},  $
that 
$$ \int_{[0,1]^n}  f\left( S_{M,1}(x) , \dots , S_{M,n}(x)
\right)
dx \ge 
\int_{[0,1]^n}  f\left( S_{M-1,1}(y) , \dots , S_{M-1,n}(y)
\right)
dy . $$
We next replace $ M $ by $M-1 $ and repeat.
Starting the process with $M= M_0$ and stopping at 
$M = 1$ yields the claimed inequality
$$ \int_{[0,1]^n}  
f\left( S_{M_0,1}(x) , \dots , S_{M_0,n}(x)\right) dx \ge 
f\left( \int_{[0,1]^n}  P(v) \right) . $$
\endproof
\paragraph{Proof of Theorem~\ref{tatarconj} :}
Choose $v \in L^p(\bR^n, \bR^n)$ and a sequence $v_r \in L^p(\bR^n, \bR^n)$ so that 
\eqref{16okt0710} and \eqref{16okt0711} hold.
Let $C_0^+ ( (0,1)^n)$ 
denote the continuous, non-negative and compactly supported functions on 
 the open unit cube $ (0,1)^n .$  We first show 
the conclusion \eqref{5dez071} under the additional restriction that 
\begin{equation}\label{5dez072}
v_{| (0,1)^n} = const, \quad \text{  and  }\quad \vp \in C_0^+ ( (0,1)^n).
\end{equation}
Clearly we may then assume that $v_{| (0,1)^n} = 0 ,$ since otherwise we replace 
$f$ by $f(\cdot + c ).$ Next we choose a smooth function 
$\a \in  C_0^+ ( (0,1)^n ) $ so that $\a(x) = 1 $ for $x \in \supp \vp . $
By considering the sequence $(\a v_r )$ instead of $ (v_r)$ we may 
further assume that  
\begin{equation}\label{16okt07z}
v_r \rightharpoonup 0 \text{  weakly in  } L^p \text{  and  } \cA_0 ( v_r ) \to 0 \text{  in  } W^{-1,p}.  
\end{equation}
By \eqref{16okt07z}
we obtain for $ v_r = ( v_{r,1},\dots, v_{r,n} ) $ that 
$$ \lim_{r\to \infty} \|R_i(v_{r,i}) \|_{L^p( \bR^n)} = 0 , \quad i\ne j .$$
Hence by \eqref{9juli3}, 
\begin{equation}\label{16okt07e}\lim_{r\to \infty}
\| v_r- P( v_r)  \|_{L^p( \bR^n, \bR^n) }
= 0 . 
\end{equation}
Since $f$ is separately convex and satisfies $ f(t) \le C ( 1 +  |t | ^p)$
we get 
\begin{equation}\label{16okt07a}
|f(s) - f(t) | \le C (1+   |s| + |t| ) ^{p-1} | s-t| . 
\end{equation}
Using \eqref{16okt07a}
and $ 1/p + 1/q = 1$ gives 
\begin{equation}\label{16okt07b}
\begin{aligned}
\int_{\bR^n}
f(v_r) \vp dx & = 
\int_{\bR^n}
f(P(v_r)) \vp dx +  \int_{\bR^n}\left(f(v_r)  -f(P(v_r)\right) \vp dx \\
&\ge \int_{\bR^n} f(P(v_r)) \vp dx - C \|1+   |v_r| + | P(v_r)| \|_p^{p/q} \| v_r - P(v_r) \|_p.
\end{aligned}
\end{equation}
Next fix $M$ and rewrite by adding and subtracting the conditional expectation operator $E_M,$
\begin{equation}\label{16okt07b}
\int_{\bR^n}
f(P(v_r)) \vp dx = 
\int_{\bR^n} 
f(P(v_r)) E_M (\vp ) dx +\int_{\bR^n} 
f(P(v_r)) ( \vp - E_M (\vp )) dx  .
\end{equation}
Clearly the conditional expectation $E_M$ satisfies 
$$\int_{\bR^n} 
f(P(v_r)) E_M (\vp ) 
dx = \int_{\bR^n} 
E_M ( f(P(v_r)))E_M (\vp )  dx . $$
Now we may invoke \eqref{18okt071}, Jensen's inequality on the range of $P.$ 
This gives,
$$\int_{\bR^n} 
E_M \left( f(P(v_r))\right) E_M (\vp ) dx  \ge \int_{\bR^n} 
f\left( E_M(P(v_r) \right) E_M (\vp )  
dx $$
Hence adding and subtracting  $f(0) $ to the leading term in the right hand side of \eqref{16okt07b}
gives   
\begin{equation}\label{16okt07d}
\int_{\bR^n} 
f(P(v_r))E_M (\vp )  dx 
\ge    \int_{\bR^n} 
f(0) E_M (\vp )  dx  +
\int_{\bR^n}  
\left( f\left(  E_M (P(v_r)\right) - f(0) \right)E_M (\vp ) 
dx
.
\end{equation}

It remains to specify  how the above estimates are to be combined:
Given $\e > 0 $ choose $M$ large enough so that 
$$ | \vp - E_M \vp | \le \e. $$
Next, depending on $M,$ and $ \e$ select $r_0 \in \bN $ so that for 
$r \ge  r_0 ,$ 
$$  | E_M ( P( v_r ) ) | \le \e \quad \text{and} \quad \|v_r - P( v_r) \|_p \le  \e . $$
Combining now  \eqref{16okt07a} -- \eqref{16okt07d} 
with our choice of $M $ and $r$ we get
$$ \int_{\bR^n}
f(v_r) \vp dx
\ge  \int_{\bR^n}f(0) \vp dx - C \e. 
$$

It remains to show how to remove the additional restriction
\eqref{5dez072}. In view of the Lipschitz condition
\eqref{16okt07a} it suffices to prove the theorem for those weak-limits
$v$ that are contained in a suitable dense set $D$ where dense refers to the 
$L^p_{\rm{loc}}$ topology.  We take 
$$D = \left\{ v \in L^p( \bR^n, \bR^n ): v \quad\text{is a finite sum of Haar functions} \right\}. $$  
Let $v \in D .$ Since the estimate \eqref{5dez071} is invariant under dilations 
$ x \to \lambda x $ it suffices to consider the case 
\begin{equation}\label{5dez075}
v(x) = \sum _{k \in \bZ^n} b_k 1_{k + (0,1)^n}(x) ,
 \end{equation}
and  only finitely many of the $b_k$ are different from zero.  

Let $\eta \in C_0^+ ( (0,1)^n)$ and extend $\eta $ to a  $(0,1)^n$
periodic continous function on $\bR^n .$ Since we proved 
\eqref{5dez071} already under the restriction
\eqref{5dez072} we obtain for functions $v$ satisfying \eqref{5dez075}
and $\vp \in  C_0^+ ( \bR^n)$ that 
\begin{equation}\label{5dez076}
\liminf_{r\to \infty} 
\int_{\bR^n} f(v_r(x)) ( \vp \cdot \eta) (x) dx
\ge \int_{\bR^n} f(v(x)) ( \vp \cdot \eta) (x) dx .
\end{equation}
Finally we remove $\eta$  from the estimate  \eqref{5dez076}.
To this end let  $\eta_k \in C_0^+ ( (0,1)^n)$ be a sequence that 
converges pointwise to $1_ {[0,1]^n}  $ and extend each $\eta_k$ periodically.
Then for each $k $ by \eqref{5dez076}
\begin{equation} 
\begin{aligned}\liminf_{r\to \infty} 
\int_{\bR^n} f(v_r(x))  \vp  (x) dx &\ge \liminf_{r\to \infty} 
\int_{\bR^n} f(v_r(x)) ( \vp \cdot \eta_k) (x) dx \\
 & \ge \int_{\bR^n} f(v(x)) ( \vp \cdot \eta_k) (x) dx .
\end{aligned}
\end{equation}
Apply now the monotone convergence theorem to conclude that 
\eqref{5dez071} holds true.

\endproof

\nocite{tat1}\nocite{tat2}\nocite{tat3}\nocite{tat4}
\nocite{mur1}\nocite{mur2}\nocite{mur3}
\nocite{daco1}\nocite{smue1}\nocite{smue2}\nocite{fonmue}
\section{Multiscale Analysis of directional Haar Projections}
In this section we outline the proof of 
Theorem~\ref{th1}.
We start by performing a multiscale analysis of $ P^{(\varepsilon)}  $
with the purpose of successively resolving the discontinuities
of the Haar system. We expand $ P^{(\varepsilon)}  $
in a series of operators, where each summand corresponds to a dyadic length scale.
Thereafter we state the estimates of Theorem~\ref{th2a} and Theorem~\ref{th2b}
that quantify 
the interplay between the resolving operators and 
the inverse of the Riesz transform $R_{i_0} . $ 
Finally we show how 
the assertions of Theorem~\ref{th1} follow.

Recall that 
$\cA = \{ \varepsilon \in \{ 0,1\}^n : \varepsilon \ne (0,\dots,0)\}.$
We decompose the projection $P^{(\varepsilon)},\varepsilon \in \cA ,$
using a smooth compactly supported approximation of unity.
To this end we choose
 $ b \in C^\infty ( \bR) , $  supported in $[-1,1] , $ so that for 
$t\in \bR ,$ 
$$ 
b(t)  = b(-t) ,\,\, 0 \le b(t) \le 4, \,\,  \Lip (b) \le 8 ,\,\, 
 {\rm and } \,\,\int_{-1}^{+1} b(t)dt = 1 . $$

Let 
$$
d(x) = b(x_1)\cdot \dots \cdot b(x_n) 
- 2^{n}  b(2x_1)\cdot \dots \cdot b(2x_n), \quad x =  ( x_1, \dots, x_n). $$
Since $b$ was chosen to be even around $0,$ we have $\int_{-1}^{+1} t b(t) dt = 0 $ hence also
\begin{equation}
\label{12okt071}
\int_\bR d(x_1,\dots, x_i,\dots , x_n) x_i dx_i = 0 ,\quad\quad (1 \le i\le n).
\end{equation} 
Let  $\Delta_\ell,$ $\ell \in \bZ$
 be the self adjoint operator defined 
by convolution as  
\begin{equation}
\label{12okt072}
\Delta_{\ell}(u) =  u*d_{\ell},  \quad\text{where}\quad
d_\ell(x) = d(2^\ell x)2^{n\ell}.
\end{equation}
For $u \in L^p( \bR^n) $ we get 
$u=\sum^\infty_{\ell=-\infty}\Delta_\ell(u).$
Convergence 
holds almost everywhere and in $L^p( \bR^n) . $
Recall that $\cS $ denotes the collection of all dyadic cubes in $\bR^n . $
Let  $j \in \bZ $ and put
\begin{equation} \label{25jan068}
\cS_j = \{Q\in \cS :|Q|=2^{-nj}\}. 
\end{equation}
Let  $\ell\in\bZ, $ $\varepsilon \in \cA ,$  define 
$T_\ell^{(\varepsilon)}$ as 

$$T_\ell^{(\varepsilon)} (u)=
\sum^\infty_{j=-\infty}
\sum_{Q  \in \cS_j}
\la u , \Delta_{j+\ell}(h^{(\varepsilon)}_{Q})\ra
h_{Q}^{(\varepsilon)}|Q|^{-1}. 
$$
Since the operators $\Delta_{j+\ell}$ are self adjoint, 
$$P^{(\varepsilon)}(u)=\sum^\infty_{\ell=-\infty}T_\ell^{(\varepsilon)} (u).$$

Let  $1 \le i_0 \le n. $ 
Recall that 
$ 
\cA_{i_0} = \{ \varepsilon \in \cA : 
 \varepsilon = (\varepsilon_1, \dots  ,\varepsilon_n) 
\quad\text{and}\quad \varepsilon_{i_0} = 1 
\}. $
Let $\e \in \cA_{i_0}.$
In Section~\ref{technical} we verify that 
 $$
T_\ell^{(\varepsilon)}  R^{-1}_{i_0}=T_\ell^{(\varepsilon)}   R_{i_0} + \sum^n_{\substack{i=1\\i \ne i_0}   }
T_\ell^{(\varepsilon)} \bE_{{i_0}}\pa{_i}R_i ,$$
where $R_i $ denotes the $i-$th Riesz transform,
 $\pa{_i}$ denotes the differentiation with respect to the 
$x_i $ variable and $\bE_{{i_0} }$ the 
 integration  with respect to the $x_{i_0}-th  $ 
 coordinate, $$
\bE_{{i_0} } (f)(x) = \int^{x_{i_0}}_{-\infty}f(x_1,\dots,s,\dots,  x_n)ds ,
\quad\quad x = ( x_1, \dots,  x_n).
$$
The following two theorems record the norm estimates for the 
operators  $T_{\ell}^{(\varepsilon)}$ and $
T_\ell^{(\varepsilon)}  R^{-1}_{i_0}$ by which we obtain  the upper bounds 
for $P^{(\varepsilon)}(u)$ stated in Theorem~\ref{th1}. 
  
\begin{theor}\label{th2a}
Let $1 < p  < \infty$ and 
$1/p + 1/q = 1 $ and $\ell \ge 0 .$ 
For $\varepsilon \in\cA$  the operator $T^{(\varepsilon)}_{\ell}$ satisfies the
norm estimates, 
\begin{equation}\label{8jan3}
\| T_{\ell}^{(\varepsilon)} \|_p \le \begin{cases}
C_p 2^{-\ell/2} &\quad\text{if} \quad  p \ge 2 ;\\
C_p 2^{-\ell/q} &\quad\text{if} \quad  p \le 2 .\\
\end{cases}
\end{equation}
Let $1 \le i_0 \le n, $ and $\varepsilon \in \cA_{i_0}$ then  
\begin{equation}\label{81jan3}
\| T_{\ell} ^{(\varepsilon)}R_{i_0}^{-1} \|_p   \le \begin{cases}
C_p 2^{+\ell/2} &\quad\text{if} \quad  p \ge 2 ;\\
C_p 2^{+\ell/p} &\quad\text{if} \quad  p \le 2 .\\
\end{cases}
\end{equation}
\end{theor}

\begin{theor}\label{th2b}
Let $1 < p < \infty .$  
Let $\ell \le 0 .$ Then for $\varepsilon \in \cA $
the operator $T_{\ell} ^{(\varepsilon)}$ satisfies the
norm estimates, 
\begin{equation}\label{82jan3}
\| T_{\ell} ^{(\varepsilon)} \|_p   \le \begin{cases}
C_p 2^{-|\ell |/p} &\quad\text{if} \quad  p \ge 2 ;\\
C_p 2^{-|\ell |} &\quad\text{if} \quad  p \le 2 .\\
\end{cases}
\end{equation}
If  moreover $1 \le i_0 \le n, $ and $ \varepsilon \in \cA_{i_0} ,$
then  
\begin{equation}\label{82jan3xx}
\| T_{\ell} ^{(\varepsilon)} R_{i_0}^{-1} \|_p \le \begin{cases}
C_p 2^{-|\ell |/p} &\quad\text{if} \quad  p \ge 2 ;\\
C_p 2^{-|\ell |} &\quad\text{if} \quad  p \le 2 .\\
\end{cases}
\end{equation}

\end{theor}
We show how Theorem~\ref{th2a} and  Theorem~\ref{th2b} yield the 
proof of Theorem \ref{th1}.
\paragraph{Proof of Theorem \ref{th1}.}
Let  $1 \le i_0 \le n . $
Define $M \in \bN$  by  the relation 
\begin{equation}
\label{17feb3}
2^{M-1} \le \frac{||u||_p ||R_{i_0}||_p}{||R_{i_0} (u)||_p}\le 2^{M} .
\end{equation}
Consider first $p\ge 2. $ Let $\varepsilon \in \cA_{i_0}.$ 
Theorem~\ref{th2a} and  Theorem~\ref{th2b}
imply that
$$\sum^\infty_{\ell=M}|| T_{\ell} ^{(\varepsilon)} ||_p \le C_p  2^{-M/2}
\quad\text{and}\quad
\sum^{M-1}_{\ell= -\infty}|| T_{\ell} ^{(\varepsilon)} R_{i_0}^{-1}
||_p \le C_p 2^{M/2}.
$$
Since 
$
P ^{(\varepsilon)}(u)=\sum^\infty_{\ell=-\infty}T_\ell ^{(\varepsilon)} (u) 
$
 triangle inequality gives that 
\begin{equation}
\begin{aligned}
 \|P ^{(\varepsilon)}(u)\|_p &\le  \sum ^{\infty}_{\ell=M}||T_\ell ^{(\varepsilon)}\|_p\|
  u||_p+\sum^{M-1}_{\ell= -\infty}||T_\ell ^{(\varepsilon)} R^{-1}_{i_0}||_p\, ||R_{i_0}(u)||_p\\
&\le   C_p    2^{-M/2}||u||_p +  C_p  2^{M/2}||R_{i_0}(u)||_p.
\end{aligned}
\end{equation}
 Inserting the value of $M$ specified in \eqref{17feb3} gives
$$  C_p 2^{-M/2}||u||_p +  C_p  2^{M/2}||R_{i_0}(u)||_p \le  
C_p ||u||_p^{1/2}||R_ {i_0}(u)||_p^{1-1/2}.$$

Assume next that  $p\le 2. $ Let $q$ be the H\"older conjugate index 
to $p$ so that $1/p+ 1/q = 1 . $
By Theorem~\ref{th2a} and  Theorem~\ref{th2b}, for 
$\varepsilon \in \cA_{i_0},$ 
$$\sum^\infty_{\ell=M}|| T_{\ell} ^{(\varepsilon)} ||_p \le C_p  2^{-M/q}
\quad\text{and}\quad
\sum^{M-1}_{\ell= -\infty}|| T_{\ell} ^{(\varepsilon)} R_{i_0}^{-1}
||_p \le C_p 2^{M/p}.
$$
Triangle inequality applied to
 $
P^{(\varepsilon)}(u)=\sum^\infty_{\ell=-\infty}T_\ell^{(\varepsilon)} u 
$
gives
\begin{equation}
\begin{aligned}
 \|P ^{(\varepsilon)}(u)\|_p &\le  \sum ^{\infty}_{\ell=M}||T_\ell ^{(\varepsilon)}\|_p\|
  u||_p+\sum^{M-1}_{\ell= -\infty}||T_\ell ^{(\varepsilon)} R^{-1}_{i_0}||_p\, ||R_{i_0}(u)||_p\\
&\le   C_p    2^{-M/q}||u||_p +  C_p  2^{M/p}||R_{i_0}(u)||_p.
\end{aligned}
\end{equation}
With $M$ defined as in \eqref{17feb3} above we obtain 
$$ C_p    2^{-M/q}||u||_p +  C_p  2^{M/p}||R_{i_0}u||_p
\le    C_p 
     ||u||_p^{1/p}||R_{i_0}u||_p^{1-1/p}.
$$

\endproof

\section{Tooling up}
\label{technical}
In this section we prepare the tools provided by 
the Calderon Zygmund School of Harmonic Analysis. 
They simplify our tasks and save the reader time and effort.
We exploit the   Haar system indexed by (and supported on)
dyadic cubes,  its unconditionality in $L^p ( 1 < p < \infty ), $
projections onto block bases of the Haar system, 
the connection of  singular integral operators
to wavelet systems, and interpolation theorems for  operators on dyadic $H^1$
and  dyadic $\BMO .$ 
\paragraph{The Haar system  in $\bR^n.$}
We base this review on the work of T. Figiel \cite{figsingular}
 and Z. Ciesielski \cite{zC87}.
Denote by $\cD $ the collection of all dyadic interval in
the real line $\bR ,  $  and let 
$\{h_I :  I \in \cD \} $ be the associated $L^\infty$ normalized Haar system.
It forms a complete orthogonal system in $L^2(\bR).$
Analogs of the Haar system in the multi-dimensional case
were developed by   Z. Ciesielski in \cite{zC87}.
For our purposes the mere tensor products of the one dimensional Haar
system is not quite sufficient. Instead we employ the Haar system supported 
on
dyadic cubes.

Recall that $ \cS $ denotes the collection of dyadic cubes in 
$\bR^n . $
and that  
$ \cA = \{ \varepsilon \in \{ 0, 1 \}^n : 
\varepsilon \ne ( 0,\dots . 0 ) \} . 
$ The system
$$ \{ h_Q^{(\varepsilon)} :  Q\in \cS, \varepsilon \in \cA \}$$
is a complete orthogonal system in $L^2(\bR^n)$ with $\| h_Q^{(\varepsilon)}\|_2^2 = |Q| . $
It is also an unconditional basis in  $L^p(\bR^n)$ $(1 < p < \infty). $ 
Given 
$f \in L^p(\bR^n)$ define its dyadic square function 
$\bS(f) $ as
 \begin{equation}
\label{3maerz2}
\bS^2(f) = 
\sum_{\varepsilon\in \cA,\,Q\in \cS}
 \la f,   h_Q^{(\varepsilon)} \ra^2  1_Q |Q|^{-2}
\end{equation}
The norm of $f \in L^p(\bR^n)$ and that of its square function $\bS(f) $
are related by the estimate
\begin{equation}
\label{4maerz3}
 C_p^{-1} \|f \|_{L^p(\bR^n)}  \le \|  \bS(f)\|_{L^p(\bR^n)} \le 
C_p \|f \|_{L^p(\bR^n)} ,
\end{equation}
where $C_p \le Cp^2/(p-1). $
Repeatedly we exploit the unconditionality of the Haar system 
in the following form. 
Let $ \{ c_Q^{(\varepsilon)} :  Q\in \cS, \varepsilon\in \cA \}$
be a bounded set of coefficients and
$f \in L^p(\bR^n).$ 
Then  
$$g = \sum_{\varepsilon\in \cA ,\,Q\in \cS}
 c_Q^{(\varepsilon)} \la f,   h_Q^{(\varepsilon)} \ra 
h_Q^{(\varepsilon)}|Q|^{-1}
, $$
satisfies the square function estimate $\bS(g) \le  
\left(\sup| c_Q^{(\varepsilon)}|\right)\bS(f) , $
hence by  \eqref{4maerz3}
\begin{equation}
\label{4maerz4}
 \|g\|_{ L^p(\bR^n)} \le  C_p \left(\sup| c_Q^{(\varepsilon)}| \right) \cdot \|f\|_ {L^p(\bR^n)}.
\end{equation}

\paragraph{Wavelet systems.}
We refer to Y. Meyer and R. Coifman \cite{coifme} for the 
unconditionality  of the wavelet systems and the fact that 
they are equivalent to the Haar system. 
Recall that $\cS$ denotes the collection of 
dyadic cubes in $\bR^n .$
We say that
$$\{ \psi_Q^{(\varepsilon)} : Q \in \cS, \varepsilon \in \cA \} $$
 is a wavelet system  
if 
$\{ \psi_Q^{(\varepsilon)} /\sqrt{| Q|} : Q \in \cS, \varepsilon \in \cA \} $
is an orthonormal basis in $L^2(\bR^n)$
satisfying $\int \psi_Q^{(\varepsilon)} = 0 $ and
there exists $C > 0 $ so that for $ Q \in \cS,$ and $ \varepsilon \in \cA  $
the following structure condition holds,
\begin{equation}
\label{55maerz1}
 \supp\psi_Q^{(\varepsilon)} \sbe C\cdot Q,
\quad \quad |\psi_Q^{(\varepsilon)}| \le C ,
\quad\quad \Lip(\psi_Q^{(\varepsilon)}) \le C \diam(Q)^{-1}.
\end{equation}
The wavelet system 
$\{ \psi_Q^{(\varepsilon)}  : Q \in \cS, \varepsilon \in \cA \} $ is  an unconditional 
basis in $ L^p (\bR^n)$ $(1 <p < \infty)$ and 
equivalent
to the Haar system $\{ h_Q^{(\varepsilon)} : Q \in \cS, \varepsilon \in \cA \} :$
Indeed there exists  $C_p \le C  p^2/(p-1) , $
so that for any choice of finite sums,
$$
f =\sum_{\varepsilon \in \cA,\, 
Q \in \cS}   a_Q^{(\varepsilon)} h_Q^{(\varepsilon)} 
\quad\text{and}\quad g =
\sum_{\varepsilon \in \cA , \,
Q \in \cS} a_Q^{(\varepsilon)} \psi_Q^{(\varepsilon)} , 
$$ 
the following norm estimates hold, 
\begin{equation}
\label{5maerz1}
C_p^{-1} \left\| f
\right\|_{ L^p (\bR^n)} \le  
\left\| g
\right\|_ {L^p (\bR^n)} \le 
 C_p
\left\|f
 \right \|_ {L^p (\bR^n)}. 
\end{equation}

\paragraph{Notational convention.}
Given a dyadic cube $Q \in \cS $ we write
$h_Q$ as shorthand for any
of the  functions  
\begin{equation}
h_Q^{(\varepsilon)}, \, \varepsilon \in \cA .
\end{equation}
If a statement in  this paper involves $h_Q$ where
$Q\in \cS $  then that statement is meant to hold
true
with   $h_Q$
replaced by any of the functions 
$h_Q^{(\varepsilon)},\varepsilon \in \cA .
$

\paragraph{Square function estimates and integral operators.}
In this (and the following) paragraph we isolate a class of integral operators 
for which boundedness in  $L^p  (\bR^n)$ $(1<p<\infty)$ can be obtained 
directly from the 
 unconditionality of the Haar system.
(Naturally we discuss those operators here because they will appear 
in later sections.)
Let $\{ c_Q^{} , Q \in \cS \} $ be a set of bounded coefficients where  
(for convenience) only  finitely many of them are $\ne 0 . $
Let $u \in  L^p (\bR^n).$
Then
\begin{equation}
\label{5maerz8}  
K(u) (x) = \int_{\bR^n} k(x,y) u(y) dy \quad\text{with kernel}\quad
k(x,y) = \sum_{ Q \in \cS}  c_Q^{}  h_Q^{}(x)  h_Q^{}(y) |Q|^{-1}, 
\end{equation}
satisfies the square function estimate $ \bS (K(u)) \le 
\left(\sup | c_Q^{}| \right) \bS (u) .$ Hence by
 \eqref{4maerz4},
\begin{equation}
\label{4maerz11} 
\| K(u)\|_ {L^p (\bR^n)} \le C_p\left( \sup | c_Q^{}| \right) \cdot \|u\|_ {L^p (\bR^n),} 
\end{equation}
where $C_p \le C\max\{p^2 , p/(p-1)\} . $
\paragraph{Projections onto block bases.}
Our  reference to projections  onto block bases of the Haar system
is \cite{jpw1} by P. W. Jones.
Let $\cB $ be a collection of dyadic cubes. For  $Q \in \cB$
let $\cU(Q)$ denote a collection of pairwise disjoint
dyadic cubes. We assume that the collections 
$\cU(Q)$ are disjoint as $ Q $ ranges over the cubes in $\cB .$
More precisely we assume  the 
following conditions throughout:
\begin{equation}\label{5maerz10}
\text{
If $ W \in   \cU(Q) ,\,W' \in   \cU(Q'), $
 and $Q\ne Q'$ then $W \ne W' . $}
\end{equation}
\begin{equation}\label{5maerz11}
\text{
If   $W, W' \in \cU(Q)$ and  $W \ne W'$  then $W\cap W' = \es . $}
\end{equation}

Consider the block bases 
$$
 d_Q^{} = \sum _{W \in \cU ( Q) } h_W^{},\quad Q \in \cB.
$$
Given scalars $c_Q$ we are interested in the operator
\begin{equation}
\label{5april1}
 K_1(u)  = \sum_{Q\in\cB }
c_Q \la u , h_Q^{}\ra  d_Q^{} |Q|^{-1}
\end{equation}  
that maps $\sum_{Q\in\cB }
a_Q   h_Q^{} $ to $\sum_{Q\in\cB }
a_Qc_Q   d_Q^{} .$
Similarly, given a 
wavelet system $\{ \psi_K\} $
as above and scalars $b_W$
we consider the block bases 
$$
\tilde\psi_Q = \sum_{W \in \cU (Q)}   b_W \psi_W $$
and the operator 
$$K_2(u) = \sum_{Q\in \cB } 
c_Q^{}     \left\la  u , h_Q
\right\ra \tilde\psi_Q   |Q|^{-1} . $$
We shall see below that $K_2$ can be controlled by $K_1. $
To estimate $K_1 (u)$ it is sometimes convenient to use 
a different collection of cubes as follows. 
Let $U(Q)= \bigcup_{Q \in \cU(Q)} W$ denote the pointset
covered by the collection $ \cU(Q).$ 
Suppose that there exist dyadic cubes
$E_1(Q), \dots, E_k(Q)$ , where $k$ may depend on $Q,$
so that 
 $$U(Q)
\sbe E_1(Q)\cup  \dots \cup  E_k(Q).
$$ 
Assume
that  the collections 
$\{E_1(Q), \dots, E_k(Q)\}$
are  disjoint as $ Q$ ranges over the cubes in 
$ \cB .$ 
Let 
\begin{equation}\label{5maerz12}
g_Q^{}=
\sum_{i = 1 }^{k} 
h_{ E_i(Q)}^{},\quad Q \in \cB, 
\end{equation}
put  $\g = \sup | c_Q^{}| , $
and define the  integral operator
$$
 K_0(u) =   \g 
 \sum_{Q\in\cB }
\la u , h_Q^{}\ra  g_Q^{} |Q|^{-1} .
$$
Our construction gives  the square function estimate
$$\bS(K_1 (u)) \le \bS(K_0(u) ),$$ 
hence 
$\|K_1 (u) \| _ p \le C_p  \|K_0(u)\|_p.$
Consequently, $L^p-L^q$ duality gives the norm estimate
 \begin{equation}
\label{5april2}
 \|K_1 ^* \| _ p \le C_p  \|K_0^*\|_p .
\end{equation}
Note that the transposed operators $K_1^*$ and $K_0^*$ are given as,  
$$
K_1^*(u)= \sum_{Q\in \cB} c_Q^{}\la u , d_Q^{}  \ra h_Q^{}|Q|^{-1}
\quad\text{and}\quad
K_0^*(u) =  \g \sum_{Q\in \cB} \la u , g_Q^{}  \ra h_Q^{}|Q|^{-1}.
$$ 
\paragraph{Exchanging Haar functions and  wavelets.}
The equivalence of the wavelet system to the Haar basis 
allows us to write down further examples of $L^p $ bounded integral operators.
We use again the notational convention to 
write
$\psi_Q$ denoting
any of the  wavelet functions  
$
\psi_Q^{(\varepsilon)}, \varepsilon \in \cA .
$

Assume that $ \cU(Q) ,$ $ Q \in \cB $ satisfies 
\eqref{5maerz10}
and \eqref{5maerz11}.
Let  
$ b_W , W \in \cU (Q) $
be   scalars, and assume that  $|b_W|\le B . $ 
Recall that
$$K_2(u) = \sum_{Q\in \cB } 
c_Q^{}     \left\la  u , h_Q
\right\ra \tilde\psi_Q   |Q|^{-1},  \quad\quad 
\tilde\psi_Q = \sum_{W \in \cU (Q)}   b_W \psi_W ,$$
and that  $K_1$ was  defined in  \eqref{5april1}. 
Since $K_2$ can be viewed as the composition of $K_1$ with the map
$h_W \to b_W \psi_W$ it follows 
from \eqref{4maerz4} and \eqref{5maerz1}   that 
\begin{equation}
\label{5april3}
\left\|
K_2(u)
\right\|_ {L^p (\bR^n)} \le C_p B
\cdot \|K_1(u)\|_{L^p (\bR^n)}. 
\end{equation}
Duality gives estimates for the transposed operator as, 
\begin{equation}
\label{5april4}
\left
\|
K_2^*\right\|_p  \le 
C_p B
\|K_1^*\|_p ,  
\end{equation}
where 
\begin{equation}
\label{5april5}
K_2^*(u) = \sum_{Q\in \cB } 
c_Q^{}   \la u , \tilde \psi_Q 
\ra h_Q|Q|^{-1}
\quad\text{and}\quad
K_1^*(u)= \sum_{Q\in \cB} c_Q^{}\la u , d_Q^{}  \ra h_Q^{}|Q|^{-1}.
\end{equation}

\paragraph{Calderon Zygmund kernels.}
We use the book by   Y. Meyer and R. Coifman \cite{coifme} 
as our source for singular integral operators
and their  relation  to wavelet systems. 
Let $\{ k_Q^{} : Q \in \cS  \} $
be a family of functions satisfying $\int k_Q^{} = 0  $
and these standard estimates: There exists $C>0$ so that for $Q \in \cS ,$
\begin{equation}
\label{6april11}
 \supp k_Q^{} \sbe C\cdot Q, \quad\quad | k_Q^{}|\le 1 ,
\quad\quad \Lip(k_Q^{}) \le C \diam(Q)^{-1}.
\end{equation}
Let $\{ c_Q^{} : Q \in \cS  \} $
be a bounded sequence of scalars. Assume for simplicity that 
only finitely many of the $c_Q^{}$ are different from zero. 
Then
$$
k_3(x,y) = \sum  c_Q^{}  \psi_Q^{}(x)  k_Q^{}(y) |Q|^{-1}, $$  
defines a  standard Calderon-Zygmund kernel (see \cite{coifme})
so that 
$$ K_3(u) (x) = \int k_3(x,y) u(y) dy$$
satisfies the norm estimate
$$\|  K_3(u)\|_p \le C_p \sup|c_Q^{}|\cdot \|u \|_p. $$
By 
\eqref{5maerz1}, the operator
$$ K_4(u) (x) = \int k_4(x,y) u(y) dy \quad\text{with kernel}\quad
k_4(x,y) = \sum  c_Q^{}  h_Q^{}(x)  k_Q^{}(y) |Q|^{-1}, $$
satisfies 
\begin{equation}
\label{6april13}
  \|K_4(u)\|_p\le C_p \sup|c_Q^{}|\cdot \|u \|_p.
\end{equation}

We will apply \eqref{6april13} in the following specialized situation.
Let $W$ be a dyadic cube and let $V$ be a cube in $\bR^n$ (not necessarily dyadic) 
so that 
\begin{equation}
\label{5sept072}
 V \spe C_1\cdot W ,\quad \quad|V| \le C_2 |W|. 
\end{equation}
Let $Q\sbe W$ be a dyadic cube. Since 
$\int k_Q = 0$ and $\supp k_Q \sbe V $ 
we have
$$
\la u, k_Q \ra = \la 1_V ( u - m_V (u)) , k_Q\ra , $$
where $  m_V (u) = \int_V u   / |V| . $
This yields the identity
$$
\sum _{Q\sbe W} \la u, k_Q \ra h_Q |Q|^{-1} =
\sum _{Q\sbe W}\la 1_V ( u - m_V (u)) , k_Q\ra h_Q |Q|^{-1}. $$
To the kernel 
$ \sum _{Q\sbe W}  h_Q^{}(x)  k_Q^{}(y) |Q|^{-1}$
we apply the estimate \eqref{6april13} with $p = 2 . $
Since the Haar system is orthogonal we obtain
\begin{equation}
\label{4sep071}
\begin{aligned}
\sum _{Q\sbe W}\la u, k_Q \ra^2|Q|^{-1}
&= \|\sum _{Q\sbe W} \la u, k_Q \ra h_Q |Q|^{-1}\|_2^2 \\
&= \|\sum _{Q\sbe W} \la 1_V ( u - m_V (u)) h_Q \ra |Q|^{-1}\|_2^2\\
&\le \|  1_V ( u - m_V (u)) \|_2^2.
\end{aligned}
\end{equation}
With \eqref{4sep071}  we obtain $\BMO$ estimates for operators with  Calderon Zygmund  kernels as above.
\paragraph{The Riesz Transforms.}
We review  basic facts about  Riesz transforms and base 
the discussion 
on chapter III of \cite{eS70} by E. M. Stein. 
Let $\cF$ denote the Fourier transformation on $\bR^n . $
The  Riesz transform 
$R_i$ is a Fourier multiplier defined by    
\begin{equation}\label{25jan069}
\cF(R_i(u))(\xi) = -\sqrt{-1} \frac{\xi_i}{|\xi|}\cF(u)(\xi)
\quad\text{where}\quad 1 \le i\le n , \quad \xi = ( \xi_1,\dots , \xi_n ).
\end{equation}
Riesz transforms, satisfy the estimates  $ \|R_j u \|_p \le C_p \| u \|_p  $
$(1 < p <\infty),$ 
hence define  bounded linear operators on the reflexive $L^p( \bR^n)$ spaces.
The defining relation \eqref{25jan069} yields a convenient 
formula for the inverse of $R_i,$ again  by Fouriermultipliers.
Consider for simplicity $i = 1. $
Let $u$ be a smooth and compactly supported  test function
such that  $\cF^{-1} ( |\xi|/\xi_1 \cF(u)(\xi)) $ is well defined. 
Then compute $\cF(R_1^{-1}(u))(\xi)$
as 
$$
\begin{aligned}\cF(R_1^{-1}(u))(\xi)&=-\sqrt{-  1}
\cF(u)(\xi)\frac{|\xi|}{\xi_1}=\cF(u)(\xi)\frac
{-\sqrt{-  1}}{\xi_1}\sum^n_{i=1}\frac{\xi_i^2}{|\xi|}\\
&=-\sqrt{-  1}\cF(u)(\xi)\left[\frac{\xi_1}{|\xi|}  +  
\sum^n_{i=2}
\frac{\xi_i}{\xi_1} \cdot \frac{\xi_i}{|\xi|}
\right].
\end{aligned}
$$
Taking the inverse Fourier transform 
yields
\begin{equation} \label{8jan11}
 R^{-1}_1= R_1+\sum^n_{i=2}\bE_{1}\pa{_i}R_i ,
\end{equation}
where
$\bE_{1} (f)(x_1,\dots , x_n) = \int^{x_1}_{-\infty}f(s,x_2,\dots , x_n)ds$
and $\pa{_i}$ denotes the partial differentiation with respect to the 
$i -th$ coordinate.  

Next fix $ 1 \le i_0 \le n $ and $\varepsilon \in \cA_{i_0} . $ 
After permuting the coordinates  the above  calculation gives 
 the formula for $R^{-1}_{i_0}$ as follows
\begin{equation}\label{18sept062}
 R^{-1}_{i_0}= R_{i_0} + \sum^n_{\substack{i=1\\i \ne i_0}   }
\bE_{{i_0}}\pa{_i}R_i . 
\end{equation}

\paragraph{Dyadic $\BMO ,$ $H^1_d$ and  Interpolation.}
We use  \cite{BS88} by C. Bennett and R. Sharply  as  basic reference to 
interpolation theorems. 
Recall first the definition of dyadic $\BMO.$ 
Let $f \in L^2(\bR^n)$ 
with Haar expansion given by \eqref{10juli1}
We say that $f$ belongs to dyadic $\BMO  $
and write $ f \in \BMO_d$
if the norm defined by \eqref{3maerz3} is
finite
 \begin{equation}
\label{3maerz3}
\|f\|_{\BMO_d}^2
=  \left|\int f\right| ^2 + \sup_{Q\in \cS} \frac{1}{|Q|}   \sum_{\varepsilon \in \cA}
\sum_{ W \sbe Q } \left\la f,   h_W^{(\varepsilon)} \right\ra^2 |W|^{-1}.
\end{equation}
Given a dyadic cube $Q$ the system
$$ \{ 1_Q\} \cup \{  h_W^{(\varepsilon)}:\, W \in \cS, W\sbe Q , \varepsilon \in \cA \} $$
is a complete orthogonal system in the Hilbert space $L^2 ( Q , dt ) . $
This yields the identity
$$ 1_Q ( f - m_Q(f)) = \sum_{\varepsilon \in \cA}\sum_{W\sbe Q} 
\la f,   h_W^{(\varepsilon)} \ra  h_W^{(\varepsilon)} |W|^{-1} , $$
where $ m_Q(f) = (\int_Q f)/| Q| . $ Hence the $\BMO _d $ norm of $f$ can be rewritten as 
\begin{equation}
\label{6sep071}
\|f\|_{\BMO_d}^2
=  
 \left| \int f \right| ^2 + 
\sup_Q \int_Q | f(t) - m_Q(f)|^2\frac{dt}{|Q|} .
\end{equation}

Given $f 
\in \BMO _d$
with $\int f = 0 . $ Let 
$\cG = \{ W \in \cS : \, \exists  \varepsilon \,  \,\la f,   h_W^{(\varepsilon)}\ra  \ne 0\}.$ It is well known that 
in order to evaluate   the $\BMO _d $ norm of $f$ it suffices to consider the cubes in $\cG . $
Put
$$
\begin{aligned}
A_0 & 
& = \sup_{Q\in \cG}\frac{1}{|Q|}\sum_{\varepsilon \in \cA}
\sum_{ W \sbe Q } 
\la f,   h_W^{(\varepsilon)} \ra
^2 |W|^{-1}.
\end{aligned}
$$
We claim that \begin{equation}\label{4sept073}
 A_0 = \|f\|_{\BMO_d}^2 . 
\end{equation}
It suffices to observe that 
 $ A_0 \ge  \|f\|_{\BMO_d}^2 ,  $ since  $ A_0 \le  \|f\|_{\BMO_d}^2 ,  $ by definition.
To this end we fix a dyadic cube $K\in \cS$ so that 
$K \not\in \cG . $ Let $\cM \sbe \cG $ denote the collection of 
maximal cubes of $\cG $ that are contained in $K. $
(Maximality is with respect to inclusion.) Thus $\cM$
consists of pairwise disjoint dyadic cubes,
$$\sum_{Q \in \cM } |Q| \le |K|, $$
and, 
$$ \sum_{\varepsilon \in \cA}\sum_{ W \sbe K } 
\la f,   h_W^{(\varepsilon)} \ra
^2 |W|^{-1}
= \sum_{Q \in \cM }\sum_{\varepsilon \in \cA}\sum_{ W \sbe Q } 
\la f,   h_W^{(\varepsilon)} \ra
^2 |W|^{-1}.
$$
Since $\cM \sbe \cG ,$ for $Q \in \cM,$
$$ \sum_{\varepsilon \in \cA}\sum_{ W \sbe Q } 
\la f,   h_W^{(\varepsilon)} \ra
^2 |W|^{-1} \le A_0|Q| . $$
Consequently we have the following estimates
$$
\begin{aligned}
 \sum_{\varepsilon \in \cA}\sum_{ W \sbe K } 
\la f,   h_W^{(\varepsilon)} \ra
^2 |W|^{-1}
& = A_0 \sum_{Q \in \cM } |Q| \\
& =  A_0 |K|.
\end{aligned}
$$
Taking the supremum over all such $K$  implies that   $ A_0 \ge  \|f\|_{\BMO_d}^2 .  $

We review the definition of dyadic $H^1,$  its
relation to the scale of $L^p $ spaces and to $\BMO_d.$ 
Let $K$ be a dyadic cube  in $\bR^n. $ 
We say that $a : \bR^n \to \bR $ is a dyadic atom if
\begin{equation}
\label{5maerz15}
\|a\|_ {L^2(\bR^n)} \le |K|^{-1/2},\quad \supp a \sbe K,
\quad\text{and} \quad\int a = 0. 
\end{equation}
By definition a function 
$f \in L^1 (\bR^n)$ belongs to dyadic $H^1 $ if there exists a sequence of 
dyadic atoms 
$\{a_i\}$ and a
sequence of scalars $\{\lambda _i \}$
so that 
\begin{equation}
\label{55maerz15}
 f = \sum \lambda _i a_i \quad\text{and}\quad \sum |\lambda _i| < \infty.
\end{equation}
We denote
\begin{equation}
\label{5maerz16}
\|f\|_{H^1_d} = \inf\{ \sum |  \lambda _i| \}
\end{equation}
where the infimum is extended over all  representations 
\eqref{55maerz15}.
For the resulting space of functions we write $H^1_d.$
Recall also that the dual Banach space to $H^1_d$ is identifiable with
 $\BMO_d.$

Interpolation of operators links
  the spaces $H^1_d,$ $\BMO_d$ on the one hand  and 
the scale of $L^p $ spaces on the other hand. 
Assume that $T$ is a bounded operator on $H^1_d$
and on $L^2. $ Let $A_1$ denote the the norm of 
 $T$ on $H^1_d$ and let $A_2 $ denote the norm of 
 $T$ on $L^2 .$ Then for $ 1 < p <2 $ and $\theta = 2- 2/p $
$$ \|T\|_p \le CA_1^{1-\theta} A_2^\theta. $$
If on the other hand the operator $T$ is bounded on  $\BMO_d$
with norm equal to $A_\infty $ then  for $ 2 < p <\infty $ 
and $\theta =  2/p $ 
$$ \|T\|_p \le CA_\infty^{1-\theta} A_2^\theta. $$

In addition to dyadic $\BMO $ at one point of the proof we employ
the continuous analog of $\BMO_d . $
Let $f \in L^2(\bR^n).$ Let $W\sbe \bR^n $ be a 
cube (not necessarily dyadic). Write
$$ m_W(f) = \int_W f(t) \frac{dt}{|W|}. $$
We say that $f \in \BMO(\bR^n) $ if
$$\|f\|_{\BMO(\bR^n)}^2
= \left| \int f \right| ^2 + \sup_W \int_W | f(t) - m_W(f)|^2\frac{dt}{|W|} < \infty ,
$$
where the supremum is extended over all
cubes   $W\sbe \bR^n $ 
(not just dyadic ones).
Clearly for a given function
$\|f\|_{\BMO(\bR^n)} \ge \|f\|_{\BMO_d} .$
In Section~\ref{basic}
we use  $\BMO(\bR^n)$ and interpolation as follows.
Let $T: L^2(\bR^n) \to  L^2(\bR^n)$
and $T: \BMO_d \to \BMO(\bR^n) $ be bounded.
Let $A_2$ be the operator norm of $T$ on 
 $L^2(\bR^n)$ and 
put 
$$ A_\infty = \|T: \BMO(\bR^n) \to \BMO_d  \|. $$
Then for $1 < p < \infty  $ and $\theta = 2/p , $
$$  \|T\|_p \le CA_\infty^{1-\theta} A_2^\theta. $$

\section{Basic Dyadic Operations}
\label{basic}
The norm estimates for the operators 
$T_\ell^{(\varepsilon)}$
 reflect boundedness of  two basic dyadic operations. 
These are  rearrangement operators of the Haar basis
and averaging  projections onto block bases of the Haar system.
In this section we isolate the basic dyadic models 
and prove estimates in the spaces $ H^1 ,$ $L^2 $ and $\BMO . $
In later sections the boundedness properties of  
$T_\ell^{(\varepsilon)} , \ell \le 0 ,$ 
are reduced to the case of rearrangement operators. 
 The estimates for $T_\ell^{(\varepsilon)} , \ell \ge 0 ,$
are harder and involve rearrangements as well as 
orthogonal projections onto certain ring domains,
surrounding the discontinuity set of Haar functions.
\subsection{Projections and Ring Domains}\label{ringdomains}

The following definitions enter in the construction
of the orthogonal projection \eqref{25jan0611}. 
Recall the set of directions 
$\cA = 
\{ \varepsilon \in \{ 0, 1 \}^n : \varepsilon \ne ( 0,\dots . 0 ) \}.$
Let $\cB$ be a collection of dyadic cubes. For $ Q \in \cB $
and $\varepsilon\in \cA$ let $ D^{(\varepsilon)}(Q) $ 
denote the set of discontinuities of the Haar function
$h_Q^{(\varepsilon)}. $ Fix 
$\lambda \in \bN $ and define 
$$ D^{(\varepsilon)}_\lambda(Q) = \{ x \in \bR^n : \dist 
( x ,D^{(\varepsilon)}(Q))
\le C 2^{-\lambda}\diam (Q)  \} . $$
Thus  $D^{(\varepsilon)}_\lambda(Q) $ is the set of points that have 
distance $\le C 2^{-\lambda}\diam (Q)$
 to the set of discontinuities of  $h_Q^{(\varepsilon)}. $
Let $k(Q) \le  C2^{\lambda (n-1)} $ and let 
$E_1(Q), \dots, E_{k(Q)}(Q)$ be the  collection of all dyadic cubes 
satisfying 
\begin{equation}
\label{11juli3}
 \diam(E_k(Q)) =  2^{-\lambda}\diam (Q), \quad\quad E_k(Q)\cap
 D^{(\varepsilon)}_\lambda(Q) 
\ne \es . 
\end{equation}
We {\it assume} throughout this chapter that 
$\cB $ is such that the  collections $\{E_1(Q), \dots, E_{k(Q)}(Q) \} $
are pairwise disjoint as $ Q$ ranges over $\cB . $

Thus we defined a covering of $ D^{(\varepsilon)}_\lambda(Q)$
with dyadic cubes $\{E_1(Q), \dots, E_{k(Q)}(Q) \} $ satisfying 
these conditions:
\begin{enumerate}
\item There holds the measure estimate 
\begin{equation}\label{11juli4}
|E_1(Q) \cup \cdots \cup  E_{k(Q)}(Q)| \le C  2^{-\lambda }|Q|. 
\end{equation}
\item Let $ Q, Q_0 \in \cB , $ $k \le k(Q) $ and  $k_0  \le k(Q_0) . $
\begin{equation}\label{11juli5}
\text{
If $ E_k(Q) \sb  E_{k_0}(Q_0) $ then $ Q \sb  Q_0 . $
} 
\end{equation}
\item  Let $ Q, Q_0 \in \cB , $ $k \le k(Q) ,$  $k_0  \le k(Q_0)  $
and $ Q \sb  Q_0 . $
\begin{equation}\label{11juli6}
\text{
If $ E_k(Q) \cap  E_{k_0}(Q_0) \ne \es $ then $ E_k(Q) \sb  E_{k_0}(Q_0) .$ 
} 
\end{equation}
\end{enumerate}
Note that our hypothesis \eqref{11juli4}--\eqref{11juli6} 
are modeled after Jones's compatibility condition in \cite{jpw1}.
With 
$\cU(Q)= \{E_1(Q), \dots, E_{k(Q)}(Q) \} $ we define the block bases as 
$ g_Q = \sum_{E\in \cU (Q)} h_{E}.$
The associated projection operator is given  by the equation
\begin{equation}\label{25jan0611}
 S( u ) = \sum_{ Q \in \cB} \la u , h_Q \ra g_Q |Q|^{-1} .
\end{equation}
Recall that 
$h_Q$ is shorthand  for any
of the Haar functions 
$h_Q^{(\varepsilon)},$ 
where $\varepsilon \in \cA . $
Moreover,
if a statement in  this paper involves $h_Q$
then that statement is meant to hold
true
with   $h_Q$
replaced by any of the functions 
$h_Q^{(\varepsilon)}.$ 

The   norm estimates for the operator $S$ are recorded in the next theorem. 
For its use in the later sections of this paper
the relation between the spaces, on which the operator acts,
 and the
dependence 
of the operator norm on the value of  $\lambda $ 
becomes crucial. 
\begin{theor}
\label{projection}
There exists $C_0 = C_0( C , n ) $ so  
that the orthogonal projection given by \eqref{25jan0611}
satisfies these estimates
$$
\|S\|_{H^1_d} \le C_0 2^{-\lambda/2},\quad\quad
\|S\|_2 \le C_0 2^{-\lambda/2}, \quad\text{ and }\quad
\|S\|_{\BMO _d} \le C_0. 
$$
\end{theor}
\proof
The proof splits canonically into three parts.
The first  part  treats  $L^2 $,  
the second part   $H^1_d,$  and the last part the
$\BMO_d$ estimate of the operator $S .$
\paragraph{Part 1.} We start with  $L^2 .$
Since $|E_1(Q) \cup \cdots \cup  E_{k(Q)}(Q)| \le C_n  2^{-\lambda }|Q|, $
we have  $\|g_Q\|_2 ^2 \le   C_n  2^{-\lambda }|Q| .$
As we assume that  the  collections $\{E_1(Q), \dots, E_{k(Q)}(Q) \} $
are pairwise disjoint as $ Q$ ranges over $\cB , $
the induced block bases 
 $ \{g_Q : Q \in \cB \} $ are orthogonal.
Hence 
\begin{equation}
\label{11juli7}
\begin{aligned} 
\|S(u) \|_2^2 &= \sum_{Q \in \cB}  \la u , h_Q \ra^2 \| g_Q\|_2^2 |Q|^{-2} 
\\&\le 
C  2^{-\lambda }\|u\|_2^2.
\end{aligned}
\end{equation}

\paragraph{Part 2.} The $H^1 _d$ estimate. Let $a$ be a dyadic atom supported on a dyadic cube 
$K$ so that $\| a\|_2 ^2 \le |K|^{-1}.$
 If $ \la a , h_Q \ra \ne 0 , $ then $Q \sbe K $ and $\supp g_K \sbe C \cdot K .$
Hence 
$$ \supp S(a) \sbe C \cdot K . $$
The $L^2$ estimate \eqref{11juli7} gives 
$ \|S(a) \|_2^2  \le   C_n  2^{-\lambda }|K|^{-1} .$
As  $ \supp S(a) \sbe C \cdot K ,  $ we obtain 
the 
  $H^1 _d $ estimate, $ \|S(a) \|_{H^1 _d} \le  2^{-\lambda /2}C . $
\paragraph{Part 3.}
The $\BMO_d$ estimate.
 Define 
$$ \cG = \bigcup_{Q \in \cB }  \{E_1(Q), \dots, E_{k(Q)}(Q) \} .$$
Given  $u \in \BMO_d ,$
 by  \eqref{4sept073},     it is sufficient to test the 
$ \BMO_d $ norm of $S(u) $ using only the cubes $ K \in \cG . $
Indeed,
$$
\|S(u) \|_{\BMO_d}^2 = \sup_{K \in  \cG} \frac{1 }{|K|} \int_K | S(u) - 
 \frac{1 }{|K|} \int_K  S(u) |^2 .
$$
Let $ K \in \cG . $ 
Note that,
$
\frac{1 }{|K|} \int_K | S(u) - 
 \frac{1 }{|K|} \int_K  S(u) |^2
$ coincides with 
\begin{equation}
\label{11juli8}
\sum_{ Q \in \cB} 
\left \la u , \frac{h_Q}{|Q|} \right\ra^2 \sum_{ \{ k : E_k(Q) \sbe K  \}
} |E_k(Q)| . 
\end{equation}
Choose $ Q_0 \in \cB , k_0 \le k( Q_0 ) $
so that 
$K = E_{k_ 0}(Q_0). $
By \eqref{11juli5}, if $ Q \in \cB $ and $  E_k(Q) \sbe  E_{k_ 0}(Q_0), $ then
$Q \sbe  Q_0 $
and if moreover $ E_k(Q) \cap  E_{k_ 0}(Q_0) \ne \es $ then, 
by \eqref{11juli6}, $  E_k(Q) \sbe  E_{k_ 0}(Q_0). $
Hence if $Q \sbe  Q_0 $ then
$$\sum_{ \{ k : E_k(Q) \sbe K
\}} |E_k(Q)|
= \int_K  g_Q^2 , $$
and 
  \eqref{11juli8} equals,
\begin{equation}
\label{11juli9}
\sum_{ Q \in \cB,\,Q \sbe Q_0} 
 \left\la u , \frac{h_Q}{|Q|} \right\ra^2  \int_K  g_Q^2 . 
\end{equation}

To get estimates for  \eqref{11juli9} consider 
 $ s \in \bN \cup \{0\}$ such that $ s \le \lambda .$
Split the (effective) index set in \eqref{11juli9} into 
$$\cH_s = \left\{ Q \in \cB :
\, 
Q \sbe Q_0 , \,\diam(Q) = 2^{-s} \diam (Q_0) ,
\,
 \int_K  g_Q^2 \ne 0 
\right\}, \quad s \le \lambda ,
$$
and
$$\cH_\infty = \left\{ Q \in \cB :
\,
 Q \sbe Q_0 ,\, \diam(Q) < 2^{-\lambda}
 \diam (Q_0) ,\,
 \int_K  g_Q^2 \ne 0 
\right\}.
$$
First estimate the contribution to \eqref{11juli9} coming from 
$\cH_\infty . $ If $Q \in  \cH_\infty $  then by \eqref{11juli4},  
 $\int_K  g_Q^2 \le C 2^{-\lambda} |Q|. $ 
Since clearly the pointset covered by  $\cH_\infty $ 
is contained in $C\cdot K , $
we get
\begin{equation}
\label{11juli11}
\begin{aligned}
\sum_{ Q \in \cH_\infty} 
 \left\la u , \frac{h_Q}{|Q|} \right\ra^2  \int_K  d_Q^2 
& \le C   2^{-\lambda} \sum_{ Q \in \cH_\infty} \left\la u , h_Q \right\ra^2  
|Q|^{-1}\\
& \le C  2^{-\lambda} \|u\|_{\BMO _d}^2 |K|.
\end{aligned}
\end{equation}
Next turn to the 
$\cH_s , s \le \lambda . $ The analysis is parallel to the previous
case.
The cardinality of $\cH_s $ is bounded by $C_n $  with $C_n$ independent
of $ s $ or $ \lambda . $
For  $Q \in \cH_s $ we get  $\int_K  g_Q^2  \le C  2^{-s}|K| . $
Hence
 $$
\sum_{ Q \in \cH_s} 
 \left\la u , \frac{h_Q}{|Q|} \right\ra^2  \int_K  g_Q^2 
 \le C  2^{-s} \|u\|_{\BMO _d}^2 |K|.
$$
Taking the sum over $  0 \le s \le \lambda , $ gives
\begin{equation}
\label{11juli10}
 \sum_{s= 0 } ^  \lambda
\sum_{ Q \in \cH_s} 
 \left\la u , \frac{h_Q}{|Q|} \right\ra^2  \int_K  g_Q^2 
 \le C   \|u\|_{\BMO _d}^2 |K|.
\end{equation}
Adding \eqref{11juli11} and  \eqref{11juli10}  
gives the $ \BMO _d $ estimate $\|S(u) \|_{ \BMO _d }
\le  C   \|u\|_{\BMO _d} . $
\subsection{Rearrangement Operators}\label{rearrangements}

We next turn to defining the rearrangement operator
$S$ given by \eqref{25jan0612} below. 
Let $\lambda \in \bN $ and let $Q\in \cS $ be a dyadic cube. 
The $\lambda -th $ dyadic predecessor of $Q,$ denoted 
$Q^{(\lambda)}, $ is given by the relation 
$$Q^{(\lambda)} \in \cS , \quad |Q^{(\lambda)}| = 2^{n\lambda}|Q|,\quad 
Q \sb Q^{(\lambda)}. $$
Let $\tau :\cS \to \cS $ be the map that associates to 
each $ Q \in \cS $ its $\lambda -th $ dyadic predecessor. 
Thus
$$\tau(Q) =  Q^{(\lambda)},\quad Q \in  \cS .$$
Clearly $\tau  :\cS \to \cS $ is not injective.
We canonically split $ \cS = \cQ_1 \cup \cdots \cup  \cQ_{2^{n \lambda}} $
such that the restriction of $\tau $ to  each of the collections
 $ \cQ_k ,$
is injective:
Given $ Q \in \cS ,$ form
$$ \cU(Q)  = \left\{ W \in \cS : W^{(\lambda)} = Q \right\} . $$
Thus $ \cU(Q) $ is a covering of $Q$ and 
contains exactly  $2^{n \lambda}$ pairwise disjoint dyadic cubes.
We enumerate them,
rather arbitrarily, as
$ W_1(Q) , \dots ,  W_{2^{n \lambda} }(Q). $
For $ 1 \le k \le  2^{n \lambda} ,$  define 
$$ \cQ_k = \left\{ W_k(Q) : Q \in \cS \right\} . $$
Note that  $\tau :\cQ_k \to \cS $ is a bijection, and 
$$ \tau (  W_k(Q) ) = Q , \quad\quad  W_k(Q)\in \cQ_k,
\quad Q \in \cS  . $$

Let  $ 1 \le k \le  2^{n \lambda} .$ Let 
$\{ \vp_Q^{(k)} : Q \in \cS \} $ be a family of functions 
for which $ \int \vp_Q^{(k)} = 0$ and which satisfy the following 
structural conditions: There exists $C>0 $ so that for each $Q \in \cS $
\begin{equation}\label{11juli12}
 \supp \vp_Q^{(k)} \sbe C \cdot Q,\quad\quad | \vp_Q^{(k)} | \le C , 
\quad\quad \Lip (  \vp_Q^{(k)} ) \le C \diam (Q) ^{-1} . 
\end{equation}
We emphasize that the  actual function $  \vp_Q^{(k)}$  may  depend
on $k ,$ by contrast 
  the structural conditions \eqref{11juli12} 
 are independent of  the value of
$k . $ Define the operator $S$ by the equation
\begin{equation}
\label{25jan0612}
 S(g) = \sum_{k = 1 }^{ 2^{n \lambda} }
\sum_{Q \in  \cQ_k } \left\la g ,  \vp_{\tau(Q)}^{(k)}\right\ra h_Q |Q|^{-1}.
\end{equation}
The action of $S$ is best understood by 
viewing it as the transposition of the rearrangement operator 
defined by $\tau $ followed by a  Calderon Zygmund Integral.
The next theorem records the 
 operator norm of $ S,$ particularly  its joint $(n,\lambda) -$dependence, 
on the spaces $H^1_d ,$ $L^2$ and $\BMO_d . $
\begin{theor}
\label{11julitheorem1}
The operator  $ S$ defined by \eqref{25jan0612} is bounded 
on the spaces $H^1_d ,$ $L^2$ and from $\BMO (\bR^n)$ to 
 $\BMO_d . $ The 
norm estimates depend on the value of $\lambda \in \bN$ and the dimension of the 
ambient space $\bR^n $ as follows:
\begin{equation}\label{11juli20}
 \|S\|_2 \le  C_0  2^{n \lambda},\quad\quad
 \|S\|_{H^1_d} \le  C_0  2^{n \lambda},\quad
\quad
  \|S: \,\BMO (\bR^n) \to  \BMO _d \,\| \le C_0  \lambda ^{1/2}  2^{n \lambda}.
\end{equation}
\end{theor}
\proof
The three parts of the proof correspond to the 
three operator estimates in \eqref{11juli20}. The first
part treats $L^2 ,$ the second part $H^1_d$ and the third part $\BMO_d.$ 
\paragraph{Part 1.}
We start with $L^2 $. 
Let $ u \in L^ 2 . $
Then 
$$ \|S(u)\|_2^2 = \sum_{k = 1 }^{ 2^{n \lambda} }
\sum_{Q \in  \cQ_k } \left\la u ,  \vp_{\tau(Q)}^{(k)}\right\ra^2  |Q|^{-1}.
$$
Let $1\le k \le 2^{n\lambda}.$
Since 
$ \tau : \cQ_k \to \cS $ is bijective, 
the standard conditions \eqref{11juli12} and the $L^2 $ estimates for 
Calderon Zygmund operators \eqref{6april13} yield,
\begin{equation}\label{11juli13}
\sum_{Q \in  \cQ_k } \left\la u ,  \vp_{\tau(Q)}^{(k)}\right\ra^2 
 |\tau(Q)|^{-1}
\le C \|u\|_2^2 .
\end{equation}
Recall that   $|\tau(Q)| = 2^{n\lambda}|Q| . $
On the left hand side of \eqref{11juli13} 
 replace  $|\tau(Q)|^{-1}$ by $2^{-n\lambda}|Q|^{-1}
$ then  take the sum  over  
$1\le k \le 2^{n\lambda}.$ This gives 
$$\sum_{k = 1 }^{ 2^{n \lambda} }
\sum_{Q \in  \cQ_k } \left\la u ,  \vp_{\tau(Q)}^{(k)}\right\ra^2  |Q|^{-1}
\le C 2^{2n\lambda}\|u\|_2^2.$$ 
Hence  $\|S\|_2 \le  C_0  2^{n \lambda},$ as claimed.
\paragraph{Part 2.}
The $H^1_d$ estimate. Let $a$ be a dyadic atom supported on a dyadic cube $K.$
Define
$$
\cH = \left\{ Q \in \cS :\, \diam (\tau(Q)) \ge \diam (K) ,\, 
\la a ,  \vp_{\tau(Q)}^{(k)}\ra \ne 0 \right\}.
$$
Then put $S(a) = b_1 + b_2 $ where 
$$b_1 =
\sum_{Q \in  \cH  } \left\la S(a) ,  h_Q 
\right\ra h_Q |Q|^{-1},
$$
and $ b_2 =  S(a) - b_1 . $
We treat  separately the norm of $b_1 $ and $b_2 . $
First we estimate $\|b_1\|_{H^1_d}. $ Fix $s \in \bN  \cup \{0\} $
and put
$$
\cH_s = \left\{ Q \in \cH :\, \diam (\tau(Q)) =  2^{s}\diam (K) \right\} .
$$
Let $ Q \in \cQ_k \cap \cH_s  $ and let $q \in Q . $
As $\int a = 0  $ we obtain
$$
\begin{aligned}
\left|\left\la a , \vp_{\tau(Q)}^{(k)}\right\ra\right| & = \left|\left\la a , 
\vp_{\tau(Q)}^{(k)} - \vp_{\tau(Q)}^{(k)}(q)\right\ra\right|\\
& \le C \|a\|_{ L^{1}} \diam (Q) \Lip ( \vp_{\tau(Q)}^{(k)})
\end{aligned}
$$
By the structural conditions \eqref{11juli12},  $ Q \in \cQ_k \cap \cH_s  $
implies  $  \Lip ( \vp_{\tau(Q)}^{(k)}) \le C 2^{-s} \diam (K)^{-1}. $
Hence $|\la a , \vp_{\tau(Q)}^{(k)}\ra| \le  C 2^{-s} . $
Note that the cardinality of $\cQ_k \cap \cH_s $ is 
bounded by an absolute constant $C.$
Hence,
\begin{equation}\label{11juli14}
\sum_{s = 0 }^{\infty}\sum_{k = 1 }^{ 2^{n \lambda} }
\sum_{Q \in  \cQ_k  \cap \cH_s}
| \la a ,  \vp_{\tau(Q)}^{(k)}\ra|
\le C 2^{n\lambda}.
\end{equation}
Since $ h_Q/|Q| $ is of  norm one in $ H^1_d ,$ the
triangle inequality and \eqref{11juli14} give
$ \|b_1\|_{ H^1_d} \le C 2^{n\lambda}.$ 
It remains to consider $ \|b_2\|_{ H^1_d}.$
Here the estimates are   a direct consequence of the 
operator $L^2 $ norm of  $S.$
First 
$$
\begin{aligned}
\|b_2\|_2^2 & \le \|S(a)\|_2\\
            & \le C 2^{2n\lambda}\|a\|_2^2\\
            & \le C  2^{2n\lambda} |K|.
\end{aligned}
$$
Second, a moments reflection shows that the 
Haar support of $b_2$ is contained in 
$C \cdot K . $ 
Let $$\cM = \{W \in \cS : W \cap \supp b_2  \ne 0 , \,\, |W| = |K| \}$$
Clearly the union of the cubes in $\cM $ covers  $\supp b_2.$
The cardinality of $\cM $ is bounded by a constant $C_n ,$
and $\int _{W} b_2 = 0 $ for $ W \in \cM . $
 Hence  the functions 
$$  C^{-1}  2^{-n\lambda}1_{W}  b_2, \quad W \in \cM , $$
are dyadic atoms, 
and 
$ \|b_2\|_{ H^1_d} \le C 2^{n\lambda}.$
Since  $\|S(a)\|_{ H^1_d} \le \|b_1\|_{ H^1_d} + \|b_2\|_{ H^1_d}$
it follows that   $\|S\|_{ H^1_d} \le  C_0 2^{n\lambda}.$
\paragraph{Part 3.} Let $ u \in \BMO( \bR^n). $
We obtain the  $\BMO _d $ estimate for $S(u)$  
by verifying that 
for every dyadic cube  $W,$
\begin{equation}\label{11juli15}
\sum_{k = 1 }^{ 2^{n \lambda} }
\sum_{Q \in  \cQ_k , Q\sbe W } \la u ,  \vp_{\tau(Q)}^{(k)}\ra^2  |Q|^{-1}
\le   C |W|\cdot  \lambda \cdot  2^{2 n \lambda}\cdot
 \|u\|_{\BMO( \bR^n)}^2.
\end{equation}
To this end
fix a dyadic cube $W.$  Split $\{Q \in  \cS , Q\sbe W \} = \cG \cup \cH ,$
where
$$ \cH = \{ Q \in \cS :\, Q \sbe W,\, \diam(Q) \ge \diam (W) 2^{-\lambda} \} 
\quad\text{and}\quad  \cG = \{Q \in  \cS , Q\sbe W \} \sm  \cH . $$
Fix 
$1 \le k \le  2^{n \lambda},$   put $ \cG_k =  \cG \cap \cQ_k  $
and 
observe that $$\bigcup_{Q \in \cG_k} \tau(Q) \sbe W . $$
Recall further that $\tau :  \cG_k \to \cS $ is injective. 
Hence the standard conditions \eqref{11juli12}, 
the  Calderon-Zygmund estimate \eqref{4sep071}, 
and \eqref{5sept072}
yield
\begin{equation}\label{11juli16}
 \sum_{Q \in  \cG_k  } \la u ,  \vp_{\tau(Q)}^{(k)}\ra^2  |\tau(Q)|^{-1}
\le   C  |W|\cdot   \|u\|_{\BMO( \bR^n)}^2.
\end{equation}
Next  replace $|\tau(Q)|^{-1}$ by $2^{-n\lambda}|Q|^{-1} ,$ 
then  take the sum  of \eqref{11juli16} over  
$1\le k \le 2^{n\lambda} .$ 
We obtain that 
$$ \sum_{k=1}^{2^{n\lambda}}\sum_{Q \in  \cG_k  } 
\la u ,  \vp_{\tau(Q)}^{(k)}\ra^2  |Q|^{-1}
\le   C  |W|\cdot  2^{2 n \lambda}\cdot
 \|u\|_{\BMO( \bR^n)}^2.$$

We turn to estimating the contribution
to \eqref{11juli15} coming from $\cH . $ Let $ 0 \le s \le \lambda . $
Write
$$ \cH _s = \{ Q \in \cH : \diam ( Q) = 2^{-s} \diam ( W )\} .$$ 
The cardinality of $ \cH _s$ equals $ 2^{ns} . $ 
It is useful to observe that,  since  $s \le \lambda , $
there exists exactly one dyadic cube $K_s $ so that 
$$ \tau ( Q) = K_s , \quad \text{for all} \quad Q\in \cH_s . $$
Hence the following  identity holds
\begin{equation}\label{11juli17} 
\sum_{k=1}^{ 2^{n\lambda}}\sum_{Q \in  \cH_s \cap \cQ_k  } 
\la u ,  \vp_{\tau(Q)}^{(k)}\ra^2  |Q|^{-1}
=
\la u ,  \vp_{K_s}^{(k)}\ra^2  
\left[
\sum_{Q \in   \cH_s} |Q|^{-1}
\right]
. 
\end{equation}
Each  $Q \in   \cH_s $ satisfies  $ |Q| = |W| 2^{-ns} .$
As  $ \cH _s$ has cardinality equal to  $ 2^{ns} , $
it follows that 
$$\sum_{Q \in   \cH_s} |Q|^{-1} = 2^{2ns} |W|^{-1}. $$
By definition
$ |K_s| =  2^{-ns + n\lambda} |W| .$ Squaring 
and 
regrouping gives
$$  2^{2ns} |W|^{-1} = 2^{2 n\lambda} |K_s|^{-2} |W| . $$
Hence the right hand side of \eqref{11juli17} equals
\begin{equation}\label{11juli18}
 2^{2 n\lambda}  |W| 
\left\la u ,  \vp_{K_s}^{(k)} \right \ra^2  |K_s|^{-2} .
\end{equation}
By \eqref{11juli12}, $\|\vp_{K_s}^{(k)} \|_2 \le  |K_s| ^{1/2} . $
 Let $B_s $ be  a cube in $\bR^n $ so that 
$\supp ( \vp_{K_s}^{(k)}) \sbe B_s $ and $\diam (B_s) \le C \diam (K_s) . $
Let $ m_{B_s} (u) = \frac{1}{|B_s|} \int_{B_s}u(x)dx . $
 As $\int \vp_{K_s}^{(k)} = 0 $ we get 
\begin{equation}
\label{28juli1}
\begin{aligned}
|\la u ,   \vp_{K_s}^{(k)} \ra | 
& = | \la u  - m_{B_s}( u ) , \vp_{K_s}^{(k)} \ra | \\
& \le C\|1_{B_s} \cdot ( u -   m_{B_s} (u) ) \|_2   |K_s| ^{1/2}\\
& \le  C |K_s| \cdot  \|u\|_{\BMO( \bR^n)}.
\end{aligned}
\end{equation}
Inserting  \eqref{28juli1} into 
  \eqref{11juli18}  gives  that the latter is  bounded by 
\begin{equation}\label{11juli19}
 C 2^{2 n\lambda}  |W| \cdot \|u\|_{\BMO( \bR^n)}^2.
\end{equation}

Thus we showed that the left  hand side of 
\eqref{11juli17} equals 
\eqref{11juli18} which in turn is bounded by  \eqref{11juli19}.
Hence
\begin{equation}
\label{25jan0614}
\sum_{k=1}^{ 2^{n\lambda}}\sum_{Q \in  \cH_s \cap \cQ_k  } 
\la u ,  \vp_{\tau(Q)}^{(k)}\ra^2  |Q|^{-1}\le C2^{2 n\lambda}  |W| \cdot \|u\|_{\BMO( \bR^n)}^2.
\end{equation}
Finally in \eqref{25jan0614} we
take the sum over 
 $ 0 \le s \le \lambda  $ and obtain  \eqref{11juli15}
\endproof

\section{The Proof of Theorem~\ref{th2a}.}
In this section we prove Theorem~\ref{th2a}.
The sub-sections ~\ref{sub1} -- \ref{sub3} are devoted to the estimates for the operator
$T_\ell^{(\varepsilon)},$ $\ell\ge 0.$ 
In sub-section ~\ref{sub4} we discuss the reduction of the 
estimates for $T_\ell^{(\varepsilon)} R^{-1}_{i_0},$ $ \varepsilon 
\in \cA_{i_0},$ to those of 
$T_\ell^{(\varepsilon)}.$
Recall that 
$$ 
\cA_{i_0} = \{ \varepsilon \in \cA : 
 \varepsilon = (\varepsilon_1, \dots  \varepsilon_n) 
\quad\text{and}\quad \varepsilon_{i_0} = 1 
\}. $$
Let $\varepsilon \in \cA_{i_0}.$
Let $\ell\ge 0 .$
Recall that   
 for $j \in \bZ $  we let   $\cS_j $ 
 be  the collection of all dyadic cubes  in 
$\bR^{n}$  with measure equal to $2^{-nj}.$ 
Let $Q\in \cS_j $ and define
\begin{equation}
\label{28mai1}
 f^{(\varepsilon)}_{Q,\ell} = \Delta _{j+\ell} ( h_Q^{(\varepsilon)}). 
\end{equation}
With the abbreviation \eqref{28mai1} we have
\begin{equation}
\label{27mai7}
T_\ell^{(\varepsilon)} (f) = \sum_{Q \in \cS } \la f ,f^{(\varepsilon)}_{Q,\ell}\ra 
 h_Q^{(\varepsilon)} |Q|^{-1} .
\end{equation}
The functions  $f^{(\varepsilon)}_{Q,\ell}$
have  vanishing mean  and satisfy the basic estimates
\begin{equation}
\label{27mai6}
\supp  f^{(\varepsilon)}_{Q,\ell} \sbe D_\ell^{(\varepsilon)}(Q),\quad
| f^{(\varepsilon)}_{Q,\ell}| \le C, \quad \Lip ( f^{(\varepsilon)}_{Q,\ell})
\le C2^\ell (\diam (Q))^{-1},
\end{equation}
where $D^{(\varepsilon)}_\ell(Q) $ is the set of points that have 
distance $\le C 2^{-\ell}\diam (Q)$
 to the set of discontinuities of  $h_Q^{(\varepsilon)}. $
{\em Based only} on the expansion \eqref{27mai7}
and the scale invariant  conditions  \eqref{27mai6} we prove in the
following subsections that $T_\ell^{(\varepsilon)} ,\,\ell \ge 0 $ satisfies the norm estimates
\begin{equation}
\label{28mai3}
\| T_{\ell}^{(\varepsilon)} \|_p \le \begin{cases}
 C_p 2^{-\ell/2} \quad\text{for}\quad p \ge 2;\\
 C_p 2^{-\ell /q} \quad\text{for}\quad p \le 2 .
\end{cases}
\end{equation}

To this end we 
 decompose the 
operator $T_\ell^{(\varepsilon)},$ $\ell\ge 0$ into a series of operators 
$T_{\ell, m}, m \in \bZ$  using 
 a
wavelet system 
$\{ \psi_K^{(\a)} : K \in \cS, \a \in \cA \} $
so that $\{ \psi_K^{(\a)} /\sqrt{| K|}
\} $
is an orthonormal basis in $L^2(\bR^n),$
satisfying  $\int \psi_K^{(\a)} = 0 $ and
the structure conditions,
$$
 \supp\psi_K^{(\a)} \sbe C\cdot K,
\quad \quad |\psi_K^{(\a)}| \le C ,
\quad\quad \Lip(\psi_K^{(\a)}) \le C \diam(K)^{-1}.
$$
To simplify expressions below
 we suppress the superindeces  $(\a)$   and, with a slight 
abuse of notation,  
in place  of  $\{\psi_K^{(\a)}\}$ we write 
just $\{ \psi_K \}.$ Then  
 expanding a function $f $ along the wavelet basis 
we get 
$$f=\sum_{K \in \cS} \left\la f,\frac{\psi_K}{|K|}\right\ra\psi_K.$$ 
Fix $m\in\bZ$ and define  $T_{\ell,m}$ by the equation
\begin{equation}\label{8jan8}
T_{\ell,m}(f)=\sum^\infty_{j=-\infty}
\sum_{Q \in \cS_j  }
\sum_{K\in  \cS_{j +\ell+m} }
\left\la
      f,\frac{\psi_K}{|K|}\right\ra\la\Delta_{j+\ell}(h^{(\varepsilon)}_{Q}),\psi_K\ra
h_{Q}^{(\varepsilon)}|Q|^{-1}.
\end{equation}
Then 
\begin{equation}\label{8jan9}
T_{\ell}^{(\varepsilon)} (f) = \sum^\infty_{m=-\infty}T_{\ell,m} (f).
\end{equation}
In this section
we prove that  
\begin{equation}
\label{28feb1}
\sum^{-\ell - 1}_{m=-\infty}||T_{\ell,m}||_p  \le
  C_p 2^{-\ell},
\quad\text{and}\quad 
\sum_{m=-\ell}^\infty
\| T_{\ell,m} \|_p \le \begin{cases}
 C_p 2^{-\ell/2} \quad\text{for}\quad p \ge 2;\\
 C_p 2^{-\ell /q} \quad\text{for}\quad p \le 2 .
\end{cases}
\end{equation}
The bounds of  \eqref{28feb1} imply the norm estimates  for 
$T_{\ell}^{(\varepsilon)},\,\ell \ge 0 $ as stated in  
\eqref{28mai3}.   

There are three 
relevant length scales in the 
series  
\eqref{8jan8}.
\begin{enumerate}
\item The scale $2^{-j} . $ This is the sidelength of $Q \in \cS _j ,$
the cube under consideration.
\item The scale $2^{-(j+\ell)} . $ This is the scale of 
$\Delta_{j+\ell}(h^{(\varepsilon)}_{Q}).$ More precisely, since 
$\Delta_{j+\ell}$ is given by a convolution kernel of  zero mean,
the function $\Delta_{j+\ell}(h^{(\varepsilon)}_{Q})$ is supported in a strip 
 of width proportional to $2^{-(j+\ell)}  $ around the discontinuity set of 
$ h_Q^{(\varepsilon)} . $
\item The scale  $2^{-(j+\ell+m)} . $ This is the scale of the 
test functions $\psi_K ,$ $ K \in \cS_{j+\ell+m} .$
\end{enumerate} 
The estimate \eqref{28feb1} 
follows from 
Proposition~\ref{pr1}, Proposition~\ref{pr3}
 and Proposition~\ref{pr2} below which deal with
the regimes 
\begin{enumerate}
\item  
$ 2^{-(j+\ell+m)} >2^{-j}  ,$
\item
$2^{-(j+\ell+m)} < 2^{-(j+\ell)} , $
\item
$2^{-(j+\ell+m)} \in [2^{-(j+\ell)} , 2^{-j}] , $
\end{enumerate}
respectively.
Accordingly we 
treat 
separately the following three cases,
$m >0,$  
$ 0 \ge m \ge -\ell ,$ 
and $ m < -\ell . $

\subsection{Estimates for  $T_{\ell, m}$,   $\ell \ge 0$,   $m <-\ell$.}
\label{sub1}
In the case when  $m <-\ell$ and  $\ell \ge 0$ we 
 have $ 2^{-(j+\ell+m)} >2^{-j}  .$ Thus  the length scale
of the test function $\psi_K $ is larger than the scale of 
$h^{(\varepsilon)}_{Q}$ when
 $Q \in \cS _j . $

We obtain in Proposition~\ref{pr1}
 the estimates for 
 $T_{\ell, m}$   from those of the rearrangement operators treated in the previous section, and 
from the fact that the wavelet bases in $L^p ( 1 < p < \infty)$ are 
equivalent to the Haar basis.
The fruitful idea of exploiting  rearrangements of the Haar system 
in the analysis of singular integral operators originates 
in T. Figiel's work \cite{figsingular}.  
(See also \cite{MR2157745} for an exposition of 
T. Figiel's approach.) 
\begin{prop} \label{pr1}
Let $1 < p < \infty$ and $1/p + 1/q = 1.$ 
For   $\ell \ge 0$,  and   $m < -\ell$ the operator
$$
T_{\ell,m}(f)=\sum^\infty_{j=-\infty}
\sum_{Q \in \cS_j  }
\sum_{K\in  \cS_{j +\ell+m} }
\left\la
      f,\frac{\psi_K}{|K|}\right\ra\la\Delta_{j+\ell}(h^{(\varepsilon)}_{Q}),\psi_K\ra h_{Q}^{(\varepsilon)}|Q|^{-1}
$$
satisfies the norm estimate
\begin{equation}
\label{23maerz1}
\| T_{\ell , m } \|_p\le 
\begin{cases}
C_p2^{m } \sqrt{-m-\ell}     & \text{for } \quad p \ge 2; \\
C_p2^{m }   & \text{for } \quad p \le 2 .\\
\end{cases}
\end{equation}
and consequently 
$$\sum^{-\ell - 1 }_{m=-\infty}||T_{\ell,m}||_p  \le
  C_p 2^{-\ell}.$$

\end{prop}
\proof 
Fix $\ell \ge 0$ and  
$-\infty<m < -\ell.$ Let $j \in \bZ$ and fix
 a dyadic cube $Q \in \cS_j . $
Then form the collection of dyadic cubes
 
$$\cU_{\ell,m}(Q)=\{K \in \cS_{j +\ell+m} :
\la\psi_K,\Delta_{j+\ell}(h^{(\varepsilon)}_{Q})\ra\ne 0\}.$$
Clearly for $T_{\ell,m}(f)$ holds the identity
\begin{equation}
\label{17feb18}
T_{\ell,m}(f)=
\sum_{j = -\infty}^{\infty}
\sum_{Q \in \cS_j  }
\sum_{K\in  \cU_{\ell,m}(Q)   }
\left\la
      f,\frac{\psi_K}{|K|}\right\ra\la\Delta_{j+\ell}(h^{(\varepsilon)}_{Q}),\psi_K\ra h_{Q}^{(\varepsilon)}|Q|^{-1}.
\end{equation}
Observe that for 
$-\infty<m < -\ell$ the cardinality of the collection 
$\cU_{\ell,m}(Q) $ is uniformly bounded.
Next for  $K\in \cU_{\ell,m}(Q) $
we prove that
\begin{equation}
\label{17feb5}
|\la \Delta_{j+\ell}(h^{(\varepsilon)}_{Q} )  ,\psi_K\ra| \le C 2^{m}|Q | .
\end{equation}
Since  
$$\int_{\bR^n}|\Delta_{j+\ell} (h^{(\varepsilon)}_{Q} ) | dx
\le C2^{-\ell } | Q| , $$
and  since  $\Delta_{j+\ell} (h^{(\varepsilon)}_{Q} )$ 
has vanishing mean, we get for $q\in Q $ 
$$
\begin{aligned}
|\la \Delta_{j+\ell}(h^{(\varepsilon)}_{Q} )  ,\psi_K\ra|
&=|
\la \Delta_{j+\ell}(h^{(\varepsilon)}_{Q} )  ,(\psi_K  - \psi_{K}( q ) ) \ra|\\
&\le C \Lip(\psi_{K}) \diam ( Q )
 \int_{\bR^n}|\Delta_{j+\ell} (h^{(\varepsilon)}_{Q} ) | dx\\
&\le  C
\frac{\diam( Q )}{\diam(K)} 2^{-\ell} | Q |. 
\end{aligned}
$$
Next recall that  $ Q \in \cS_j $ and $ K \in \cS_{j + \ell+ m} . $
Hence $\diam( Q )= \sqrt{n} 2^{-j}  $ and  $\diam (K) =\sqrt{n} 2^{-j-m-\ell} . $ 
Inserting these values 
gives  \eqref{17feb5}.

By \eqref{4maerz11}, in combination 
with  \eqref{17feb18} and \eqref{17feb5}  we obtain   that
\begin{equation}
\label{17feb6}
\| T_{\ell, m } (f)\|_p
\le 
 C_p 2^{m} \left \| \sum_{Q \in \cS }
\sum_{K\in \cU_{\ell, m}(Q)} \left\la f,  \frac{\psi_K}{|K|}\right\ra  h_{Q} \right\|_p .
\end{equation}
Recall $K\in \cU_{\ell, m}(Q)$
 satisfies 
$|K|=|Q| 2^{n(-\ell-m)} .$
Hence
$|K|^{-1} |Q | 2^{m}= 2^{(n+1)m+n\ell}. $
Thus the right hand side of  \eqref{17feb6} 
is bounded by 
\begin{equation}
\label{17feb7}
 C_p 2^{(n+1)m + n\ell}
\left \| \sum_{Q \in \cS }
\sum_{K\in \cU_{\ell,m}(Q)} \la f,\psi_K \ra h_{Q}|Q|^{-1}\right\|_p.
\end{equation}
Given $ Q \in \cS $ 
let    $K_s(Q )$
be a cube in  $ \cU_{\ell,m}(Q). $ 
As there exist at most  $C= C_n $  cubes  in  $ \cU_{\ell,m}(Q), $ 
the expression in  
\eqref{17feb7} is bounded by 
\begin{equation}
\label{17feb8}
C_p 2^{(n+1)m + n\ell}
\max_{s\le C }\left \| \sum_{Q \in \cS }
 \la f, \psi_{K_s(Q)} \ra h_{Q}|Q|^{-1}\right\|_p.
\end{equation}

Fix $s \le C$ so that the maximum in the right hand side is assumed. 
We invoke rearrangement operators to obtain good upper bounds
for \eqref{17feb8}. 
Let $\tau : \cS \to \cS $ be the map that associates 
to $ Q \in \cS $ its $( -m-\ell ) -th $ dyadic predecessor,
denoted $Q^{( -m-\ell )} . $
Thus 
$$ \tau( Q ) = Q^{( -m-\ell )} . $$
In sub-section ~\ref{rearrangements} we defined the canonical splitting of $\cS$
as 
$$\cS = \cQ_1 \cup \dots \cup \cQ_{2^{n( -m-\ell )}},$$
so that for each fixed  $ k \le 2^{n( -m-\ell )},$
the map
$\tau : \cQ_k \to \cS $ is a bijection.
Fix now $ k \le 2^{n( -m-\ell )}$ and define the family of functions
$\{ \vp_W ^{(k)} : W \in \cS \} $ by the equations
$$   \vp_{\tau(Q)} ^{(k)} = \psi_{K_s(Q)}, \quad\quad Q \in  \cQ_k . $$
Let  $A = 2^{n(-m -\ell  )} $ and  define the rearrangement operator $ S $ by 
$$ S(f) = \sum_{k = 1}^A \sum_{Q \in \cQ_k }
\left\la f,    \vp_{\tau(Q)} ^{(k)}    \right\ra h_{Q}|Q|^{-1}.
$$

What we have obtained so far can be summarized in one line as follows
\begin{equation}\label{11juli21}
 \|T_{\ell, m }(f)\|_p \le C_p 2^{(n+1)m +n \ell  } \| S(f)\|_p . 
\end{equation}
It remains to find estimates for  $\| S(f)\|_p. $
To this end observe that  the family of functions
$\{ \vp_W ^{(k)} : W \in \cS \} $ satisfies 
the
structural conditions \eqref{11juli12}: 
There exists $C>0 $ so that for each $W \in \cS $
$$
 \supp \vp_W^{(k)} \sbe C \cdot Q,\quad\quad | \vp_W^{(k)} | \le C , 
\quad\quad \Lip (  \vp_W^{(k)} ) \le C \diam (W) ^{-1} . 
$$
Hence  Theorem~\ref{11julitheorem1} applied to 
the operator $S ,$  with $\lambda = - m - \ell , $
 gives
$$ \|S\|_p \le
\begin{cases}
C_p  2^{n ( - m - \ell ) }  \sqrt{-m-\ell}     & \text{for } \quad p \ge 2; \\
C_p  2^{n ( - m - \ell ) }   & \text{for } \quad p \le 2 .\\
\end{cases}
$$
Inserting the norm estimate for $S$ into \eqref{11juli21}
and simple arithmetic implies
 \eqref{23maerz1}.
\endproof

\subsection{Estimates for  $T_{\ell,m},$ $\ell\ge 0 ,$ $m > 0.$ }
\label{sub2}
In this subsection 
 we treat the case $m >  0$ and   $\ell\ge 0$
or equivalently $ 2^{-(j + \ell+m)} < 2^{-(j+\ell)}.$
Here the length scale of the test function $\psi_K$ is finer than the scale 
of $\Delta_{j+\ell}(h^{(\varepsilon)}_{Q}) . $
We estimate the norm 
of $T_{\ell,m}$ by reduction to the projections onto ring domains.
\begin{prop} \label{pr3}
Let  $1 < p < \infty.$ and $1/p + 1/q = 1. $
For   $m\ge 0$ and   $\ell\ge 0$,
the operator
$$
T_{\ell,m}(f)=\sum^\infty_{j=-\infty}
\sum_{Q \in \cS_j  }
\sum_{K\in  \cS_{j +\ell+m} }
\left\la
      f,\frac{\psi_K}{|K|}
\right\ra\la\Delta_{j+\ell}(h^{(\varepsilon)}_{Q}),\psi_K\ra
h_{Q}^{(\varepsilon)}|Q|^{-1}.
$$
satisfies the norm estimate
\begin{equation}
\label{18feb20}
\| T_{\ell,m} \|_p \le
\begin{cases}
 C_p2^{-m} 2^{-\ell/2} \quad\text{for}\quad p \ge 2;\\
 C_p2^{-m} 2^{-\ell /q} \quad\text{for}\quad p \le 2 .
\end{cases}
\end{equation}
\end{prop}
\proof
We divide the proof into three parts. First we rewrite 
the operator  by isolating the cubes $Q \in \cS_j $
and $K\in \cS_{j+\ell+m}$ that contribute to the series 
defining $T_{\ell,m}. $  Second we define auxiliary operators 
that dominate $T_{\ell,m}. $ These turn out to be projections onto ring domains
as considered in sub-section~\ref{ringdomains}.
Finally we invoke norm estimates 
for the resulting projections onto ring domains.
\paragraph{Part 1.}  Here we rewrite $T_{\ell,m} $
by making explicit the index set $\{K\in  \cS_{j +\ell+m} \}$
that actually contributes to the series 
defining  $T_{\ell,m}. $
Fix $Q\in \cS_j $
and define the collection of dyadic cubes 
$$\cU_{\ell,m}(Q)=\{K\in \cS_{j +\ell + m}: 
\la\Delta_{j+\ell}(h^{(\varepsilon)}_{Q}),\psi_K\ra\ne 0\}.$$
Let $U_{\ell,m}(Q)$ be the pointset that is covered by the collection 
$\cU_{\ell,m}(Q) . $ Note that  $U_{\ell,m}(Q)$
is contained in the  ring domain of points that have 
 distance  $\le C2^{- \ell - j }$ 
to the  set of discontinuities of $h^{(\varepsilon)}_{Q} . $
Thus $U_{\ell,m}(Q)$
 can
be covered by at most $C 2^{(n-1)\ell } $ dyadic cubes of diameter
$ \sqrt{n} 2^{- \ell - j } .$ We denote these  cubes (that are pairwise disjoint)  by 
$E_1, \dots , E_A $ where $A= {C2^{(n-1)\ell}}. $ 
If we wish to emphasize the dependence on
 $Q $ we write $E_k= E_k(Q ) . $ 
Thus 
$$ U_{\ell,m}(Q) \sbe  \bigcup_{k = 1 }^A
  E_k(Q ) ,
\quad 
\quad\quad
\diam (E_k(Q )) 
= \sqrt{n} 2^{- \ell - j } ,  \quad A= {C2^{(n-1)\ell}}.$$
With $\cU_{\ell,m}(Q) $ as index set we define the block bases of wavelet 
functions 
$$\widetilde \psi_{Q}=
\sum_{K\in\cU_{\ell,m}(Q)}
 \left \la \Delta_{j+\ell}(h^{(\varepsilon)}_{Q}),\psi_K
\right\ra \psi_K |K|^{-1},
$$ 
by which we rewrite the 
operator $T_{\ell,m}$ as follows, 
\begin{equation}
\label{12juli1}
T_{\ell,m}(f)= 
\sum_{Q \in \cS }
\left\la f, \widetilde\psi_{Q  } \right \ra
h^{(\varepsilon)}_{Q}|Q|^{-1}.
\end{equation} 
\paragraph{Part 2.}
Here we exploit \eqref{12juli1} and 
relate the representation $T_{\ell,m}$ to its dyadic counterpart,
the projection onto ring domains. 
To this end we start by giving  pointwise estimates for the function 
 $\widetilde\psi_{Q  }. $
Fix $K \in\cU_{\ell,m}(Q) .$  Use that
$\psi_K$ has mean zero and that $\diam (K) = \sqrt  n  2^{(-j-\ell-m)} $
to obtain, 
\begin{equation}
\label{18feb3}  
\begin{aligned}
|\la\Delta_{j+\ell}(h^{(\varepsilon)}_{Q}),
\psi_K\ra|\cdot|K|^{-1}
&\le C
\diam(K)\Lip(\Delta_{j+\ell}(h^{(\varepsilon)}_{Q}))\\
&\le C\diam(K)2^{j+\ell}\\
&= C2^{-m}.
\end{aligned}
\end{equation}

Recall that  
$$\dist (U_{\ell,m} (Q) , Q ) \le  C \cdot \diam ( Q ) \quad Q \in \cS .$$ 
Hence there exists a universal $A_0\in \bN $ 
so that for  $ j \in \bZ $ the collection 
$\cS_j $ may 
split as 
$$ \cS_j^{(1)},\dots ,  \cS_j^{(A_0)} , $$
so that for $s \le A_0 $ the sets 
$\{U_{\ell,m} (Q) :  Q\in \cS_j^{(s)} \}$
are pairwise disjoint. 
Fix  
$s \le A_0 $
and  
form the collections
$$ \cB_{s} =\bigcup _{j\in \bZ } \cS_j^{(s)}. $$
As $ s\le A_0 $
is fixed, 
the collections
$\{\cU_{\ell,m} (Q) :  Q\in \cB_{s} \}$
satisfy the conditions 
\eqref{5maerz10} and \eqref{5maerz11}.
Define
$$ d_{Q}  = \sum_{K \in \cU_{\ell,m}
 (Q)} h_K ,
$$ 
and put
$$
F_s (g) = \sum_{Q\in\cB_{s}}
\la g,  d_{Q}\ra h_Q |Q|^{-1} .
$$
By \eqref{18feb3}  and \eqref{5april4}, \eqref{5april5},
$$
\|   T_{\ell,m}  \|_p \le  C_p 2^{-m}\sum_{s = 1}^{A_0} \| F_{s}\|_p.
$$
Next  we replace the operator $F_s$ by a related one 
 that is  easier to analyze.  To this end we define
for $Q \in \cB_s ,$
$$
 g_{Q} = \sum_{ k = 1 }^{A}  h_{E_{k}(Q )}, \quad\quad A= {C2^{(n-1)\ell}},
$$
where the collection of dyadic cubes $\{ E_{1}(Q ) \dots E_{A}(Q ) \}$
are defined in part 1 of the proof. The block bases 
$\{ g_{Q}:Q \in \cB_s \} $ give rise to the operators 
$G_{s}$  defined by, 
$$G_{s}(f) =
\sum_{Q\in\cB_{s}} \la f ,   g_{Q} 
  \ra 
h_{Q} |Q|^{-1}. \quad\quad 
$$
By \eqref{5april2},
$\| F_{s}\|_p \le C_p \| G_{s}\|_p.$ Hence
\begin{equation}
\label{18feb5}  
\|   T_{\ell,m}  \|_p
\le C_p 2^{-m}\sum_{s = 1}^{A_0} \| G_{s}\|_p.
\end{equation} 
\paragraph{Part 3. }
In the last part of the proof we obtain norm estimates
for  $ T_{\ell,m}$ by recalling the bounds for the projection 
 $G_{s}^*$ obtained in Section~\ref{basic}.     Fix $ s\le A_0 ,$  let 
$$ \cB = \cB_{s}\quad\text{and}\quad G= G_{s} . $$
The transposed operator  $G^*$ is just 
$$ G^*(f) = \sum_{Q\in \cB }  \left\la f, h_{Q}\right\ra   g_{Q }|Q|^{-1}. 
$$
In  part 1 of the proof, for  $ Q \in \cB  ,$ we defined
the collections $\{ E_1(Q), \dots , E_A(Q) \} .$
They satisfy conditions 
\eqref{11juli4}--\eqref{11juli6}. 
Hence we apply  Theorem~\ref{projection} 
with $S = G^ * $ and $ \lambda = \ell . $ 
By  duality this gives 
the following three norm estimates  for  $G,$
\begin{equation}
\label{18feb6}  
 \| G \|_{H^1_d} \le  C ,\quad
\| G  \|_2 \le C2^{-\ell/2} 
\quad\text{and}\quad
 \| G  \|_{\BMO_d} \le C2^{-\ell/2} . 
\end{equation}
By interpolation and  \eqref{18feb6},
   for $1 < p < \infty $ and $1/p + 1/q = 1 $   
\begin{equation}
\label{18feb21}  
\| G \|_p \le \begin{cases}
 C_p 2^{-\ell/2} \quad\text{for}\quad p \ge 2;\\
 C_p 2^{-\ell /q} \quad\text{for}\quad p \le 2 .
\end{cases}
\end{equation}
With  \eqref{18feb21}  and \eqref{18feb5} we deduce \eqref{18feb20}.

\subsection{Estimates for  $T_{\ell,m}$,  $\ell \ge 0,$ $-\ell\le m\le 0$.}
\label{sub3}
Here we analyze the operators 
$T_{\ell,m}$,  when $\ell \ge 0,$ $-\ell\le m\le 0.$
In this case the scale of the test functions $\psi_K$ lies in between the scale of the cube $Q$ 
and that of $\Delta_{j+\ell}(h^{(\varepsilon)}_{Q}) . $ Again we 
estimate $T_{\ell,m}$ by reduction to projection operators onto ring domains,
following the pattern of the previous  sub-section.
\begin{prop}\label{pr2}
Let  $1 < p < \infty.$ and $1/p + 1/q = 1. $
Let    $\ell \ge 0$ and  $-\ell\le m\le 0$ then  the operator
$$
T_{\ell,m}(f)=\sum^\infty_{j=-\infty}
\sum_{Q\in \cS_j}
\sum_{K \in   \cS_{j+\ell+m}}
 \left\la
      f,\frac{\psi_K}{|K|}\right\ra
\la\Delta_{j+\ell}(h^{(\varepsilon)}_{Q}),\psi_K\ra h_{Q}^{(\varepsilon)}
|Q|^{-1}
$$
satisfies the norm estimate
\begin{equation}
\label{23maerz10}
\| T_{\ell,m} \|_p \le\begin{cases}
 C_p2^{m/2} 2^{-\ell/2} \quad\text{for}\quad p \ge 2;\\
 C_p2^{m/2} 2^{-\ell /q} \quad\text{for}\quad p \le 2 .
\end{cases}
\end{equation}

\end{prop}
\proof
The proof splits canonically into three parts. First we 
analyze and rewrite 
$T_{\ell,m}. $ Then  we define auxiliary operators 
that dominate $T_{\ell,m}, $ and  continue with  norm estimates 
for those  operators. As above we are led to consider projections onto
ring domains. 
\paragraph{Part 1.}
Fix $\ell \ge 0$ and  $-\ell\le m\le 0.$
Let $j \in \bZ $ and choose a 
 dyadic cube $Q \in \cS_j . $
Then form the collection of cubes
$$\cU_{\ell,m}(Q)=\{K\in \cS_{j+\ell+m} :
\la\psi_K,\Delta_{j+\ell}(h^{(\varepsilon)}_{Q})\ra\ne 0\}.$$
Observe that with the above definition of the collections
$\cU_{\ell,m}(Q)$ the following identity holds
$$
T_{\ell,m}(f)=
\sum^\infty_{j=-\infty}\sum_{Q\in \cS_j}
\sum_{K \in \cU_{\ell,m}(Q) }
 \left\la
      f,\frac{\psi_K}{|K|}\right\ra
\la\Delta_{j+\ell}(h^{(\varepsilon)}_{Q}),\psi_K\ra
h_{Q}^{(\varepsilon)}|Q|^{-1}.
$$
\paragraph{Part 2.}
Fix  $Q \in \cS_j  $ and  $K\in \cU_{\ell,m}(Q) .$
To find the auxiliary operators we prove first  that 
\begin{equation}
\label{18feb1}
\left|\left\la\Delta_{j+\ell}(h^{(1,0)}_{Q}),\psi_K\right\ra\right|
 \le C2^{ m} |K|
\end{equation}
To see this make the following observation. First  note that $ |Q| = 2^{-nj}$
and 
 $\diam ( K) =\sqrt{n} 2^{-j-m-\ell} . $ Then 
observe that $\Delta_{j+\ell}(h^{(\varepsilon)}_{Q}$ is supported in the ring domain
$D_\ell (Q) $ and 
estimate
$$
\begin{aligned}
\left|\left\la\Delta_{j+\ell}(h^{(\varepsilon)}_{Q}),\psi_K\right\ra\right|
&\le C\int _K |\Delta_{j+\ell}(h^{(\varepsilon)}_{Q})| \\
& \le C |D_\ell (Q) \cap K|\\
&\le C2^{-\ell -j }(\diam ( K ))^{n-1}\\
&\le C 2^{ m} |K|.
\end{aligned}
$$

For a cube $K\in \cU_{\ell,m}(Q)  $ its distance to 
$Q$ is bounded by the  $C \diam(Q) .$
 Hence, there exists a universal $A_0$ so that 
for  $ j \in \bZ $  the collection $\cS_j $ can be split into
$$ \cS_j^{(1)},\dots ,  \cS_j^{(A_0)} , $$
so that  the sets 
$\{U_{\ell,m} (Q) :  Q\in \cS_j^{(s)} \}$
are pairwise disjoint.
Fix $s \le A_0 $
and
 form the collections
$$ \cB_{s} =\bigcup _{j\in \bZ } \cS_j^{(s)}. $$
Note  that
$\{\cU_{\ell,m} (Q) :  Q\in \cB_{s} \}$
satisfies the conditions 
\eqref{5maerz10} and \eqref{5maerz11}.
Define
$$
F_s(f)=
2^{  m }
 \sum_{Q\in \cB_s}
\la f,  d_{Q} \ra 
h_{Q}|Q|^{-1} ,\quad\quad    d_Q =  \sum_{K\in \cU_{\ell,m}(Q) }  h_{K}           .  $$ 
  The integral  estimates
\eqref{5april5},
\eqref{5april4}
 and \eqref{18feb1} imply 
$$  \|T_{\ell,m}\|_p \le C_p \sum_{s=1}^{A_0} \|F_s\|_p .$$
\paragraph{Part 3.}
It remains  to estimate  $\|F_s\|_p . $
Notice that the collections $\cU_{\ell,m}(Q),$ $Q\in \cB_s$ 
satisfy  conditions \eqref{11juli4}--\eqref{11juli6}. 
Next apply Theorem~\ref{projection}
to $S = 2^{-m} F_s^* $ and $\lambda = \ell + m .$
By duality this yields for 
$F_s $ the norm estimates 
on 
$ L^2,$ $ H^1_d$ and $\BMO_d$
\begin{equation}
\label{18feb2}
\| F_s\|_2\le C  2^{(m-\ell)/2} ,\quad
 \| F_s\|_{ H^1_d} \le C 2^{m}  ,\quad\text{and}\quad
 \| F_s \|_{ \BMO_d} \le C 2^{(m-\ell)/2} . 
\end{equation}
By interpolation  from \eqref{18feb2} we get   for $1 < p < \infty $ and $1/p + 1/q = 1 $   
that, 
$$
\| F_s  \|_p \le \begin{cases}
 C_p    2^{(m-\ell)/2} &\quad\text{for}\quad p \ge 2;\\
 C_p 2^{m/2-\ell /q} &\quad\text{for}\quad p \le 2 .
\end{cases}
$$
\endproof 

\subsection{Estimates for $T_{\ell}^{(\varepsilon)}
R^{-1}_{i_0},$  $\ell \ge 0.$}
\label{sub4}
We give  the norm estimates for $T_{\ell}^{(\varepsilon)}R^{-1}_{i_0} , \ell \ge 0,$ $ \varepsilon \in \cA_{i_0},$ and  $1\le i_0 \le n .$ 
We do this {\em by reduction} to the estimates 
for
the operator $T_{\ell} ^{(\varepsilon)} , \ell \ge 0 .$ 
Strictly speaking we discuss 
the reduction to the proof given in the previous sub sections.
We obtain a series representing $T_\ell ^{(\varepsilon)} R^{-1}_{i_0},$
analyze the shape and form of the measures  
$\bE_{{i_0}}\pa{_i}h_Q^{( \varepsilon )}$ and describe 
how the convolution operator 
$\Delta_{j+\ell}$ acts on those measures. 
In the following analysis we also collect the information needed for the estimates  of the 
$T_{\ell} ^{(\varepsilon)} R^{-1}_{i_0}$ when $\ell \le 0.$
\paragraph{The representation of $T_\ell ^{(\varepsilon)} R^{-1}_{i_0}.$}
In Theorem~\ref{th2a} and  Theorem~\ref{th2b} we aim at estimates for
$T_\ell ^{(\varepsilon)} R^{-1}_{i_0}$ when $\varepsilon\in \cA_{i_0} . $
Hence we seek  an explicit expansion for $T_\ell ^{(\varepsilon)}
 R^{-1}_{i_0}.$
By  \eqref{18sept062} we have 
\begin{equation} \label{27mai1}
 R^{-1}_{i_0}= R_{i_0} + \sum^n_{\substack{i=1\\i \ne i_0}   }
\bE_{{i_0}}\pa{_i}R_i \quad\text{and}\quad
T_\ell^{(\varepsilon)}  R^{-1}_{i_0}=T_\ell^{(\varepsilon)}   R_{i_0} + \sum^n_{\substack{i=1\\i \ne i_0}   }
T_\ell^{(\varepsilon)} \bE_{{i_0}}\pa{_i}R_i .
\end{equation}
Let $ j \in \bZ. $ Recall 
that 
$\cS_j$ denotes the family of dyadic cubes $Q$ for which $|Q| = 2^{-nj}.$
Let $ Q \in \cS_j ,$  
$ i \ne i_0 ,$ and $ \varepsilon \in \cA_{i_0}.$ Then
form
\begin{equation}
\label{27mai2}
k_Q^{( \ell, i )} = 
\Delta_{j+\ell}\left( \bE_{{i_0}}\pa{_i}h_Q^{( \varepsilon )}\right).
\end{equation}
Thus by \eqref{27mai1}
\begin{equation}
\label{21maerz1}
T_\ell ^{(\varepsilon)} R^{-1}_{i_0} (u) = T_\ell^{(\varepsilon)}
  R_{i_0} (u) +
\sum_{Q \in \cS}  \sum^n_{\substack{i=1\\i \ne i_0}   }
\la R_i (u) , k_Q^{( \ell, i )}\ra h_Q^{( \varepsilon )} |Q|^{-1}.
\end{equation}
Given the representation \eqref{21maerz1} we 
 further analyze the functions $\{k_Q^{( \ell, i )}:\,Q \in \cS \}. $
It is only at this point of our analysis that we exploit the fact that $i_0$ and
$\varepsilon$ are related by the condition  $\varepsilon \in \cA_{i_0} .$
\paragraph{The  measures  $\bE_{{i_0}}\pa{_i}h_Q^{( \varepsilon )}$ .}
We defined $k_Q^{( \ell, i )}$ by a convolution operator  applied to  
$$\bE_{{i_0}}\pa{_i}h_Q^{( \varepsilon )}, \quad\quad i \ne i_0, \quad 
 \varepsilon \in \cA_{i_0} ,$$
where $\pa{_i}$ denotes the differentiation with respect to the 
$y_i $ variable and $\bE_{{i_0} }$ denotes 
 integration  with respect to the $x_{i_0}-th  $ 
 coordinate, $$
\bE_{{i_0} } (f)(x) = \int^{x_{i_0}}_{-\infty}f(x_1,\dots,s,\dots,  x_n)ds ,
\quad\quad x = ( x_1, \dots,  x_n).
$$
Thus, $\bE_{{i_0}}\pa{_i}h_Q^{( \varepsilon )}$ admits a  convenient factorization: 
Let  $x = ( x_1, \dots,  x_n) ,$  then
\begin{equation}
\label{12okt071a}
\bE_{{i_0}}\pa{_i}h_Q^{( \varepsilon )}(x) 
= \left[\int^{x_{i_0}}_{-\infty} h_{I_{i_0}}^{ \varepsilon_{i_0} }
(s)ds 
\right]
 \, \left[\pa{_i} h_{I_i}^{ \varepsilon_i }(x_i)\right] \,
\left[ \prod \{ h_{I_k}^{\varepsilon_k} (x_k): k \notin\{i_0,i\}\}\right].
\end{equation}
The properties of the three factors appearing in \eqref{12okt071a}
are as follows.
\begin{enumerate}
\item
As  $\varepsilon \in \cA_{i_0} ,$ we have $\varepsilon_{i_0} = 1 , $
hence the first factor in \eqref{12okt071a} 
$$ x_{i_0} \to \int^{x_{i_0}}_{-\infty} h_{I_{i_0}}^{ \varepsilon_{i_0} }
(s)ds $$
is supported in the interval $I_{i_0}.$ Furthermore it is 
bounded by $|I_{i_0}|$ and piecewise linear with nodes at 
$l(I_{i_0}), m(I_{i_0})$ and $r(I_{i_0})$ and slopes $+ 1, -1 $ or $0.$
Here we let $l(I_{i_0})$ denote the left endpoint of $I_{i_0} ,$
 and $m(I_{i_0}),$  $r(I_{i_0})$ denote its midpoint, respectively its right 
endpoint.
\item The partial derivatives $\pa{_i}$ applied to  $ h_Q^{( \varepsilon )}$
induces a  Dirac measure, at each of the discontinuities of 
$ h_{I_i}^{ \varepsilon_i } . $ The resulting formulas depend 
on the value of 
$ \varepsilon_i \in \{0,1\},$ since 
$$
\begin{aligned} 
\pa{_i} h_{I_i} &= \d_{ l(I_i)} - 2\d_{ m(I_i)} 
+ \d_{ r(I_i)} ,\\
\pa{_i} 1_{1_i} &= \d_{ l(I_i)}  - \d_{ r(I_i)}.
\end{aligned}
$$
In either case, for $\vp \in C^\infty ( \bR ) $  
the above identities yield  the estimate,
\begin{equation}
\label{12okt073} 
\left| 
\la  \pa{_i} h_{I_i}^{ \varepsilon_i }, \vp \ra 
\right|
\le 2 \sup \left\{ \frac { |\vp(s) - \vp(t)|}{|s-t|}:\, s,t, \in I \right\} |I_i| . 
\end{equation}
\item
The third factor in \eqref{12okt071a} is  the function
\begin{equation}\label{12okt074}
x \to \prod \{ h_{I_k}^{\varepsilon_k} (x_k): k \notin\{i_0,i\}\}
\end{equation}
 It is piecewise constant and assumes the values $\{ -1 , 0 , +1 \}.$
When restricted to a dyadic cube  $W$ with $\diam (W) \le \diam (Q)/2 $ the factor
\eqref{12okt074} defines a constant function. 
\end{enumerate}
As a result of the above discussion $\bE_{{i_0} }\pa{_i} h_Q^{( \varepsilon )}$
is a measure  supported on $Q$ so that for any continous function on $\bR^n,$
$$ |\la \bE_{{i_0} } \pa{_i}h_Q^{( \varepsilon )} , \vp \ra | \le |Q| \cdot \|\vp \|_\infty 
\quad\text{and}\quad 
\la \bE_{{i_0} } \pa{_i}h_Q^{( \varepsilon )} , 1 \ra = 0 . 
 $$
\paragraph{The convolution $ \Delta_{j+\ell} $ acting on  $\bE_{{i_0} }\pa{_i} h_Q^{( \varepsilon )}$.}
Recall that in \eqref{12okt072} the operator $ \Delta_{j+\ell} $ is given as convolution with  $d_{j+\ell} $
so that 
$$
\Delta_{j+\ell}\left( \bE_{{i_0}}\pa{_i}h_Q^{( \varepsilon )}\right)=
\bE_{{i_0}}\pa{_i}h_Q^{( \varepsilon )}*d_{j+\ell} , 
$$
with 
\begin{equation}
\label{12okt075}
\supp d_{j+\ell} \sbe [-C2^{-(j+\ell)} ,C2^{-(j+\ell)}] ^n , \quad |  d_{j+\ell} | \le C 2^{n(j+\ell)},
\quad 
 \Lip (d_{j+\ell}) \le C 2^{(n+1)(j+\ell)}.
\end{equation}
Moreover for $1 \le i \le n $ by \eqref{12okt071}
\begin{equation}
\label{12okt076}
\int_\bR  d_{j+\ell}(x-y) y_i dy_i = 0\quad\text{and}\quad \int_{\bR^n}  d_{j+\ell}(x-y) dy = 0, \quad x \in \bR^n. 
\end{equation}

We derive next for $k_Q^{( \ell, i )}$ its structural estimates concerning
support, Lipschitz properties and pointwise bounds. It turns out that these depend critically 
on 
the value of  $ \sign (\ell) : $
\begin{enumerate}
\item
The case  $\ell \ge 0 . $  For $ Q \in \cS $
and $\varepsilon\in \cA_{i_0}$ let 
$ D^{(\varepsilon)}(Q) $ 
denote the set of discontinuities of  the Haar function
$h_Q^{(\varepsilon)}. $ Fix 
$\ell \in \bN $ and define 
$$ D^{(\varepsilon)}_\ell(Q) = \{ x \in \bR^n : 
\dist ( x ,D^{(\varepsilon)}(Q))
\le C 2^{-\ell}\diam (Q)  \} . $$
Thus  $D^{(\varepsilon)}_\ell(Q) $ is the set of points that have 
distance $\le C 2^{-\ell}\diam (Q)$
 to the set of discontinuities of  $h_Q^{(\varepsilon)}. $

Fix $x\not\in  D^{(\varepsilon)}_\ell(Q) .$ 
As we observed in the  paragraphs following \eqref{12okt071a}
there exist $A \in \{ -1,0, 1\} $ and $a \in \bR$
so that, 
\begin{equation}
\label{18okt073}
 \bE_{{i_0}}
h_Q^{( \varepsilon )}(y) = A ( y_{i_0} - a ),  \quad\text{for}\quad y \in B(x, c2^{-(j+\ell)}).
\end{equation}
Combining now \eqref{12okt075} with 
 \eqref{12okt076} and \label{18okt073} we find
\begin{equation}
\begin{aligned}
\Delta_{j+\ell}( \bE_{{i_0}}h_Q^{( \varepsilon )}) (x)&=
\int_{\bR^n}  d_{j+\ell}(x-y) \bE_{{i_0}} h_Q^{( \varepsilon )}(y) dy \\
&= A \int_{\bR^n}  d_{j+\ell}(x-y) ( y_{i_0} - a ) dy \\
& = 0 .   
\end{aligned}
\end{equation}
Since  $ \Delta_{j+\ell} $ is a convolution operator it commutes with 
differentiation, 
and we obtain for  $ x\not\in  D^{(\varepsilon)}_\ell(Q) ,  $  
\begin{equation}
\label{12okt077}
\begin{aligned}
 \Delta_{j+\ell}( \bE_{{i_0}}\pa{_i}h_Q^{( \varepsilon )})(x) &=  
\pa{_i}\Delta_{j+\ell}( \bE_{{i_0}}h_Q^{( \varepsilon )}) (x)\\
&=  0 .
\end{aligned}\end{equation}
Combining \eqref{12okt077} with \eqref{12okt075} we obtain that 
 the functions $\{k_Q^{( \ell, i )}:\,Q \in \cS ,\, i \ne i_0,\, \ell \ge 0\} $
satisfy the structural conditions
\begin{equation}
\label{27mai3}
\supp k_Q^{( \ell, i )} \sbe D^{(\varepsilon)}_\ell(Q) ,\quad 
| k_Q^{( \ell, i )}| \le C 2^{\ell} ,\quad 
\Lip (  k_Q^{( \ell, i )}) \le C 2^{2\ell} (\diam(Q))^{-1},
\end{equation}
with $C>0$ independent of $ Q\in \cS ,\,$  $i \ne i_0,\,$ or $ \ell \ge 0 .\, $
\item
The case   $\ell \le 0 . $
In this case we use \eqref{12okt073}  and \eqref{12okt075}  to see that the family 
$\{k_Q^{( \ell, i )}:\,Q \in \cS ,\, i \ne i_0,\, \ell \le 0\} ,$
satisfies the following conditions 
\begin{equation}
\label{27mai4}
\supp k_Q^{( \ell, i )} \sbe (C 2^{|\ell|})\cdot Q ,\quad 
| k_Q^{( \ell, i )}| \le C 2^{\ell(n+1)} ,\quad 
\Lip (  k_Q^{( \ell, i )}) \le C 2^{\ell( n + 2 )} (\diam(Q))^{-1},
 \end{equation}
were again  $C>0$ is independent of $ Q\in \cS ,\,$  $i \ne i_0,\,$ or $ \ell 
\le 0 .\, $
 \end{enumerate}


\begin{prop}
\label{12juli4}
Let $1 < p < \infty. $ Let $ 1 \le i \ne i_0 \le n $ and $\varepsilon \in
\cA_{i_0}.$
For $\ell \ge 0 $ the operator $X$ defined by  
$$X(f)=
\sum_{Q \in \cS}
\la
 f,   k_Q^{( \ell, i )}
\ra h^{(\varepsilon)}_Q |Q|^{-1}, $$
satisfies the norm estimates
\begin{equation}
\label{23maerz5}
||X||_{p}\le 
 \begin{cases}
C_p2^{+\ell/2} &\quad\text{if} \quad  p \ge 2 ;\\
C_p2^{+\ell/p} &\quad\text{if} \quad  p \le 2 .\\
\end{cases}
\end{equation}
\end{prop} 

\proof
Recall the expansion \eqref{27mai7} asserting that 
$$
T_\ell^{(\varepsilon)} (f) = \sum_{Q \in \cS } 
\la f ,f^{(\varepsilon)}_{Q,\ell}\ra 
 h_Q^{(\varepsilon)} |Q|^{-1} ,
$$
where   $f^{(\varepsilon)}_{Q,\ell}$
has  vanishing mean  and satisfies the basic estimates \eqref{27mai6},
$$
\supp  f^{(\varepsilon)}_{Q,\ell} \sbe D_\ell^{(\varepsilon)}(Q),\quad\quad
| f^{(\varepsilon)}_{Q,\ell}| \le C, \quad \quad\Lip ( f^{(\varepsilon)}_{Q,\ell})
\le C2^\ell (\diam (Q))^{-1},
$$
and where $D^{(\varepsilon)}_\ell(Q) $ is the set of points that have 
distance $\le C 2^{-\ell}\diam (Q)$
 to the set of discontinuities of  $h_Q^{(\varepsilon)}. $
Using only  the scale invariant  conditions  \eqref{27mai6} we proved
that
$T_\ell^{(\varepsilon)} ,\,(\ell \ge 0 )$ satisfies the norm estimates
\eqref{28mai3}, that is, 
$$
\| T_{\ell}^{(\varepsilon)} \|_p \le \begin{cases}
 C_p 2^{-\ell/2} \quad\text{for}\quad p \ge 2;\\
 C_p 2^{-\ell /q} \quad\text{for}\quad p \le 2 .
\end{cases}
$$
Observe that by \eqref{27mai3}
the functions 
 $\{2^{-\ell} k_Q^{( \ell, i )}\}$ satisfy the very same 
structure conditions \eqref{27mai6} as   
$\{f^{(\varepsilon)}_{Q,\ell}\} .$
Hence for the norm  of the 
operator $2^{-\ell} X $ there hold the same upper bounds 
as
for  
$T_\ell^{(\varepsilon)} ,\,\ell \ge 0 .$
Consequently,
the norm of   $X$
can be estimated as 
$$ 
||X||_{p}\le  \begin{cases}
C_p  2^{\ell-\ell/2} &\quad\text{if} \quad  p \ge 2 ;\\
C_p  2^{\ell-\ell/q} &\quad\text{if} \quad  p \le 2 .\\
\end{cases}
$$
\endproof
Proposition~\ref{12juli4} in combination with \eqref{8jan3}
and \eqref{21maerz1}
implies  that for $\ell > 0 ,$
$$||T_\ell^{(\varepsilon)} R^{-1}_{i_0}||_{p}\le \begin{cases}
 C_p   2^{\ell/2} &\quad\text{if} \quad  p \ge 2 ;\\
 C_p   2^{\ell/p} &\quad\text{if} \quad  p \le 2 .\\
\end{cases}
$$

\section{The Proof of Theorem~\ref{th2b}.}
In this section we prove Theorem~\ref{th2b}.
It turns out that 
for  $\ell \le 0$ the norm estimates for 
$T_\ell ^{(\varepsilon)}R^{-1}_{i_0}$ and $T_\ell^{(\varepsilon)} $ are  much simpler than 
for 
$\ell \ge 0. $ 
Indeed for $\ell < 0 $ the scale of $Q\in \cS_j$
is finer than the scale of 
 $\Delta _{j+\ell} (h_Q^{(\varepsilon)})$ and  the discontinuities  
of the Haar function are  completely smeared out. 
We can therefore   reduce the problem 
 to estimates for rearrangement operators acting on Haar functions,
treating 
$T_\ell^{(\varepsilon)} R^{-1}_{i_0}$
 and $T_\ell^{(\varepsilon)} $
simultaneously by the same method.

Let $u$ be a smooth function with vanishing mean and compact support. 
 Let 
$ i \ne i_0 $ and  $ \varepsilon \in \cA_{i_0}.$ Then
$$
T_\ell ^{(\varepsilon)} R^{-1}_{i_0} (u) 
= T_\ell^{(\varepsilon)}  R_{i_0} (u) +
\sum_{Q \in \cS}  \sum^n_{\substack{i=1\\i \ne i_0}   }
\la R_i (u) , k_Q^{( \ell, i )}\ra h_Q^{( \varepsilon )} |Q|^{-1},
$$
where 
$$
k_Q^{( \ell, i )} = 
\Delta_{j+\ell}( \bE_{{i_0}}\pa{_i}h_Q^{( \varepsilon )}), 
\quad\quad Q \in \cS_j.
$$
Since $\ell < 0 $ 
the functions 
$\{k_Q^{( \ell, i )}:\,Q \in \cS ,\, i \ne i_0,\, \ell \le 0\} ,$
satisfy  conditions \eqref{27mai4}.
Recall also that  
$$
\label{28mai17}
T_\ell ^{(\varepsilon)}(u) = \sum_{Q \in \cS } \la u ,f^{(\varepsilon)}_{Q,\ell}\ra 
 h_Q^{(\varepsilon)} |Q|^{-1}, 
$$
where 
$$
\label{28mai18}
 f^{(\varepsilon)}_{Q,\ell} = \Delta _{j+\ell} ( h_Q^{(\varepsilon)}) ,\quad\quad Q \in \cS_j.
$$
 It is easy to see that 
also  the 
family $ \{  f^{(\varepsilon)}_{Q,\ell}: \,Q \in \cS \, ,\ell \le 0 \} $ 
satisfies 
 the same structural conditions \eqref{27mai4},
that is  
\begin{equation}
\label{27mai4a}
\supp  f^{(\varepsilon)}_{Q,\ell} \sbe (C 2^{|\ell|})\cdot Q ,\quad 
|  f^{(\varepsilon)}_{Q,\ell}| \le C 2^{-|\ell|(n+1)} ,\quad 
\Lip ( f^{(\varepsilon)}_{Q,\ell} ) \le C 2^{-|\ell|( n + 2 )} (\diam(Q))^{-1}.
\end{equation}

\begin{prop} If  $\ell \le 0 $
then 
$$
||T_\ell ^{(\varepsilon)}||_{p} + ||T_\ell^{(\varepsilon)} R^{-1}_{i_0}||_{p} 
\le \begin{cases}C_p 2^{-2|\ell|/p}   \quad\text{for} \quad p \ge 2;\\
 C_p 2^{-|\ell|}   \quad\text{for} \quad p \le 2.   
\end{cases}
$$
\end{prop} 

\nocite{MR2157745}

\proof Let $1 \le i \ne i_0 \le n . $
Let $Q \in \cS. $
 Choose signs 
$\d_{Q, i } ,\e_Q\in \{+1 ,0, - 1\}$ and form 
\begin{equation}\label{28jan061}
g_{Q,\ell} =  \left[\sum^n_{i=1,\,i \ne i_0} \d_{Q, i }  k_Q^{( \ell, i )}\right] +\e_Q  f^{(\varepsilon)}_{Q,\ell}
.
\end{equation}
We emphasize that the definition of $g_{Q,\ell}$  depends on the 
choice of signs $\d_{Q, i } ,\e_Q\in \{+1 ,0, - 1\};$
nevertheless our notation suppresses this dependence.
Note that by \eqref{27mai4} and \eqref{27mai4a} the functions $\{g_{Q,\ell}\}$
are of mean zero and  satisfy  structure conditions,
not depending on the choice of signs, namely  
\begin{equation}
\label{6april20}
\supp g_{Q,\ell}
\sbe C 2^{|\ell|}\cdot Q ,\quad
|g_{Q,\ell}|
\le C2^{-(n+1)|\ell|},\quad
\Lip (g_{Q,\ell})
\le C2^{-(n+2)|\ell|}  \diam(Q)^{-1}.
\end{equation}

Consider the rearrangement 
$\tau :\cS \to \cS $ that maps $ Q\in \cS $ 
to its $|\ell| -th $ dyadic predecessor.
Let $ \cQ_1 , \dots ,\cQ_{2^{n|\ell|}}$ be the canonical splitting
of $ \cS$ so that for fixed 
$ k \le   2^{n|\ell|} $ the map 
$\tau :\cQ_k \to \cS $
is bijective. Fix $ k \le   2^{n|\ell|} .$
Determine the family 
$\{ \vp_W ^{(k)} : W \in \cS \} $ by the equations
\begin{equation}\label{28jan062} 
  \vp_{\tau(Q)} ^{(k)} = 2^{(n+1) |\ell |} g_{Q,\ell}, \quad\quad Q \in  \cQ_k . \end{equation}
Thus defined the functions $\vp_W ^{(k)} $ are of mean zero and satisfy the structural 
conditions
$$ \supp \vp_W^{(k)} \sbe C \cdot W,\quad\quad | \vp_W^{(k)} | \le C , 
\quad\quad \Lip (  \vp_W^{(k)} ) \le C \diam (W) ^{-1} . $$ 
Define the operator
$$
S(u) = \sum_{k = 1 } ^{2^{n |\ell |}} \sum_{Q \in \cQ_k } 
\left\la u , \vp_{\tau(Q)} ^{(k)}\right\ra
 h_Q^{(\varepsilon)} |Q|^{-1} .
$$
Apply  Theorem~\ref{11julitheorem1} to $S$ with $\lambda = | \ell | . $
This yields 
\begin{equation}
\label{25jan0617}
 \| S\|_2 \le  C_0  2^{n |\ell|},\quad\quad
 \|S\|_{H^1_d} \le  C_0  2^{n |\ell|},\quad\quad
  \|S: \,\BMO (\bR^n) \to  \BMO _d \,\| \le C_0  |\ell| ^{1/2}  2^{n |\ell|}.
\end{equation}
Note that by \eqref{28jan061} and \eqref{28jan062}  
the algebraic definition of the operator $S$  depends on the 
  choice of signs $\d_{Q, i } , \e_Q\in \{ +1,0,-1\},$  yet 
by \eqref{25jan0617} our estimates for $\|S\|_p$ are independent thereof.

 Let $g \in L^p .$ Depending 
on $g$ we choose   $\d_{Q, i } , \e_Q\in \{ +1,0,-1\},$
hence  $S ,$ so that 
\begin{equation}
\label{18feb31}
 \| T_\ell^{(\varepsilon)}(g)\|_p + \| T_\ell^{(\varepsilon)} R_{i_0}^{-1}(g)\|_p  \le C_p 2^{-(n+1)|\ell|}
C_p
\|S\|_p \|g\|_p. 
\end{equation}
Consequently, our  upper bounds for $\|T_\ell ^{(\varepsilon)}\|_p +\|T_\ell^{(\varepsilon)} R_{i_0}^{-1}\|_p$
follow  from \eqref{25jan0617}.
Indeed, by interpolation and the estimate $|\ell | ^{1/2}\le 2^{|\ell|/2} ,$ \eqref{25jan0617} and \eqref{18feb31}  imply that 
$$||T_\ell ^{(\varepsilon)}||_{p} +
||T_\ell ^{(\varepsilon)}R^{-1}_{i_0}||_{p} 
\le \begin{cases}C_p 2^{-2|\ell|/p}   \quad &\text{for} \quad p \ge 2;\\
 C_p 2^{-|\ell|}   \quad &\text{for} \quad p \le 2.   
\end{cases}
$$
\endproof

\section{Sharpness of  the exponents in Theorem~\ref{th1}.}
\label{sharp}
In this section we construct the examples showing that 
the exponents $(1/2,1/2) $ respectively $(1/p, 1/q)$
are sharp in the estimates of Theorem~\ref{th1}, 
\begin{equation}
\label{9juli9}
 ||P^{(\varepsilon)}(u)||_p \le
C_p \|u\|_p^{1/2}\|R_{i_0}(u)\|_p^{1/2},  \quad  p \ge 2 ,
\end{equation}
and
\begin{equation}
\label{9juli10}
  ||P^{(\varepsilon)}(u)||_p \le
C_p \|u\|_p^{1/p}\|R_{i_0}(u)\|_p^{1/q}, \quad  p \le 2 ,
\end{equation}
where $ 1 \le i_0 \le n $ and $\varepsilon \in \cA_{i_0} .$

When we say that we obtained sharp exponents in  Theorem~\ref{th1}
we mean the following:
Let  $\eta > 0 . $
Since  the Riesz transform is a bounded operator on 
$L^p (1 <p<\infty), $ 
 replacing in 
\eqref{9juli9} the pair of exponents 
$(1/2,1/2) $ by $(1/2 - \eta ,1/2 + \eta ) $ 
{\em would} lead 
to a statement that {\em implies} \eqref{9juli9},
hence would yield a stronger theorem. 
Our examples show, however, that 
improving the  exponents in the right hand side of \eqref{9juli9}
is impossible. (The same holds for \eqref{9juli10}.) Specifically we have this theorem: 
\begin{theor} \label{th200}
Let $1 \le i_0 \le n , $ and $\varepsilon \in \cA_{i_0} . $
Let $1 < p < \infty ,$
$1/p + 1/q = 1 .$
and $\eta >0. $ 
Then 
\begin{equation}
\label{14april1}
\sup_{u \in L^p} \frac{ ||P^{(\varepsilon)} (u)||_p}
{\|u\|_p^{1/2-\eta}\|R_{i_0}(u)\|_p^{1/2+\eta}}
= \infty
\quad\text{}\quad p\ge 2,
\end{equation}
and
 \begin{equation}
\label{14april1111}
\sup_{u \in L^p} \frac{ ||P^{(\varepsilon)} (u)||_p}
{\|u\|_p^{1/p-\eta}\|R_{i_0}(u)\|_p^{1/q+\eta}}
= \infty,
\quad\text{}\quad p\le 2.
\end{equation} 
\end{theor}
For simplicity of notation we 
 verify  Theorem~\ref{th200}
only in the case when $n = 2 . $ The passage to arbitrary $n \in \bN $ is
routine and left to the reader. 
Moreover we carry out the proof of Theorem~\ref{th200} with the following 
specification
 \begin{equation}
\label{9juli11}
n=2 ,\quad i_0 = 1,\quad \varepsilon = (1,0) .
\end{equation} 
Throughout this section we assume \eqref{9juli11} and put
$$ P = P^{(1,0)}. $$ 

We obtain Theorem~\ref{th200}
by exhibiting a sequence of test functions for which the quotient
in \eqref{14april1}
respectively \eqref{14april1111} tends to infinity.
On each test function we 
 prove lower $L^p $ bounds for 
the action of $P $
and upper $L^p $ estimates for $R_1 . $
In sub-section~\ref{building}
we define building blocks
$ s\otimes d $ and the
test functions $f_\e $ using  a procedure that resembles 
that of adding independent 
copies of the basic building blocks.
The   proof of \eqref{14april1}
requires  upper estimates for
$\|f_\e \|_p $ and $\|R_1 ( f_\e )\|_p ,$
that we  prove in sub-section~\ref{upper} 
and  a lower estimates for $\|P(f_\e) \|_p  $
obtained in sub-section~\ref{lower} .

\subsection{The building blocks $s\otimes d.$}
\label{building}
We build the examples showing sharpness of exponents 
on the properties of 
the functions $s\otimes d$ defined here.
Throughout this section we fix $ \e > 0 .$

Let $A,B$ be Lipschitz functions on $\bR .$ 
Assume that 
\begin{equation}
\label{4sept075}
\supp A \sbe [ 0,1] ,\quad\int A = 0 \quad {\rm and } \quad 
\supp B \sbe [-1,1] .
\end{equation}
Given $x  = (x_1, x_2) $ we define
$$
s(x_1) = A(x_1) , \quad d(x_2) = B(x_2 /\e ),  $$
$$ s\otimes d (x) = s(x_1) d(x_2). $$
We rescale $ g =  s\otimes d$  to a dyadic square  $ Q = I \times J $ 
as follows.
Let $l_I , l_J $ denote the left endpoint of 
$I$ respectively $J . $
Put
$$
s_I(x_1) = s(\frac{ x_1 - l_I}{|I|}), 
\quad
d_J(x_2) = d(\frac{ x_2 - l_J}{|J|}), 
$$ 
and
\begin{equation}\label{15juli1}
 g_Q(x) = s_I(x_1) d_J(x_2). 
\end{equation}

We next define  the testing  function  $f_{\e}$ 
that  is obtained by first 
forming ``almost independent'' copies of $g = s \otimes d $ 
and then adding $\frac{1}{\e}$ of those.
Below 
we define a collection of dyadic squares $\cG$ and form
\begin{equation}
\label{21april12}
f_\e = \sum_{ Q \in  \cG } g_Q.
\end{equation}
To define $\cG$ we proceed as follows.
\label{defining}
Fix $ j \in \bN. $ Let $\cD_j $ denote the collection of dyadic intervals 
$I$ satisfying
$$ I \sbe [0,1] \quad \text{and} \quad  |I| = 2^{-j}.$$
Let $\cL_j \sbe \cD_j $ satisfy
\begin{equation}
\label{21april1000}
I, J \in  \cL_j  \quad \text{implies } \quad \dist( I , J ) \ge |I| ,
\end{equation}
and 
\begin{equation}
\label{21april100}\sum_{ J \in  \cL_j } |J| = \frac12 . 
\end{equation}
To define  $ \cL_j $ simply take  the even numbered intervals 
of $\cD_j , $ counting from left to right.
Next assume that $ \e >0 $ is  power of $1/2,$ thus 
\begin{equation}
\label{21april11}
 \e = 2^{-n_0}  \quad \text{for some } \quad   n_0 \in \bN. 
\end{equation}
For $ 1 \le k \le 1/\e$ put
$$
\cG_k = \bigcup \{ I \times J : I \in \cD_{2kn_0},\,\,  J \in \cL_{2kn_0} \} 
\quad \text{and} \quad \cG =  \bigcup_{  k= 1 } ^ { 1/\e} \cG_k 
.$$
Observe that $ |Q| = \e^{4k}$ for $ Q \in \cG_k, $
and by \eqref{21april100}
\begin{equation}
\label{21april13}
 \sum_{ Q \in  \cG_k } |Q| = \frac12  
\quad \text{and} \quad
 \sum_{ Q \in  \cG } |Q| = \frac1{2\e} . 
\end{equation}

\subsection{Upper estimate for $\|f_\e\|_{p} $ and $\|R_1(f_\e)\|_{p} $.}
\label{upper}
We obtain our  $L^p $ estimates of   $f_\e $ by proving
an upper bound for its norm in the space dyadic $\BMO . $
These in turn follow from scale-invariant $L^2 $ estimates
and `` almost orthogonality'' of the functions 
$$\sum_{ Q \in  \cG_k  } g_Q, \quad\quad k \le \frac{1}{\e} . $$
\begin{prop}\label{21april14}
Let $f_\e $ be defined by \eqref{21april12}.
The support of $f_\e $ is contained in 
$ [-1,1]\times [-1,1] $ and
\begin{equation} 
\label{21april15}
 \|f_\e\|_{\BMO _d } \le C .
\end{equation}
Hence $  \|f_\e\|_{p} \le C_p .$ 
\end{prop}
\proof
Let  $Q_0 \in \cG $ and form 
$g =  \sum_{\{Q \in \cG ,\,Q \sbe Q_0\} } g_Q . $
The $\BMO_d $ inequality  \eqref{21april15}   is a consequence of 
   uniform $L^2 $ estimate
\begin{equation} 
\label{21april16}
\left\| g
\right\|_{L^2(\bR^2)}^2 \le C |Q_0| ,
\end{equation}
in combination  with  the  Lipschitz estimates,
\begin{equation} 
\label{21april17}
\sum_{\{Q \in \cG ,\,|Q| >| Q_0 |\}}
\|1_{Q_0}( g_Q - m_{ Q_0} ( g_Q)) \|_{L^2(\bR^2)} \le C \e|Q_0| ^{1/2} ,
\end{equation}
where $m_{ Q_0} ( g_Q) = | Q_0 |^{-1}\int_{ Q_0}  g_Q . $
In two separate paragraphs below we  will verify   that \eqref{21april16}
and   \eqref{21april17} hold. Before that  we show how these estimates 
yield \eqref{21april15}.
Let 
$$\cK = \{ W \in \cS : \exists \varepsilon \,\, \la f_\e , h_W^{(\varepsilon)}\ra \ne 0 \}$$
Let  $W$ be a dyadic square with $ |W| \le 1/4 ,  $  then 
$\int_W  f_\e  = 0 .$
Hence  for $ W \in \cK , $
$ \diam (W) \le 1 . $  
By \eqref{4sept073}, to estimate the $\BMO_d $ norm of $f_\e$
it suffices to test the cubes of $ \cK. $
Next we fix a dyadic square $ W \in \cK . $ 
Since  $ \diam (W) \le 1  $ we may choose   $k \in \bN_0 $
such that 
$\e ^{2(k+1)} \le \diam (W) \le \e ^{2k}. $
Define a decomposition of $\cG $ as 
$\cG=\cH_1 \cap \cH_2 \cup \cH_3 $ where 
$$
\cH_1 = \{ Q \in \cG : \diam (Q) = \e^{2k} , \, Q\cap 2\cdot W \ne \es \} ,
$$ 
$$
\cH_2 = \{ Q \in \cG : \diam (Q) \ge  \e^{2(k-1)} , \, Q\cap 2\cdot W \ne \es \} ,
$$
and
$$
\cH_3 = \{ Q \in \cG : \diam (Q) \le  \e^{2(k+1)} , \, Q\cap 2\cdot W \ne \es \} 
$$ 
Accordingly let
$$ g_j = \sum_{Q \in \cH_j } g_Q , \quad j \in \{1,2,3\}. $$
The cardinality of $\cH_1$  is bounded by C. Hence
$\|1_W g_1\|_2 \le C|W|^{1/2} .$ 
With  $ A = |W|^{-1}\int_W g_2, $ and triangle inequality
\eqref{21april17} gives
$ \int_W |  g_2 - A | ^2 \le C \e^2|W| . $
The estimate \eqref{21april16} implies $\|1_W g_3\|_2 \le C|W|^{1/2} .$ 
To see this let 
$\cM $ denote the maximal squares of $\cH_3 .$ The collection
$ \cM (\sbe \cH_3)$
consists of pairwise disjoint squares so that 
$$ \sum_{Q_0 \in \cM } |Q_0| \le C |W| .$$
Next write $ G_{Q_0} =\sum_{Q \in \cH_3, Q\sbe Q_0 } g_Q ,$ to obtain 
$$ g_3 = \sum_{Q_0 \in \cM }  G_{Q_0} 
\quad {\rm and} \quad 
\| g_3\|_2^2 = \sum_{Q_0 \in \cM } \|   G_{Q_0}      \|_2^2. $$
 Apply \eqref{21april16} to $ G_{Q_0} $ to obtain
$$
\begin{aligned}
\| g_3\|_2^2  & \le C \sum_{Q_0 \in \cM }  |Q_0| \\
              & \le  C |W|.
\end{aligned}
$$ 
Finally   $  \|1_W g_3\|_2 \le \|g_3\|_2 \le  C |W|^{1/2}. $

Moreover for $ t \in W $ there holds the identity
$$ f_\e(t) = g_1 (t) + g_2 ( t) + g_3( t ) . $$
Invoking the estimates for $g_1 ,  g_2 , g_3  $ we obtain
$$ \int_W | f_\e - A |^2 \le C |W|. $$
By  \eqref{6sep071} this estimate yields \eqref{21april15}.

\endproof

\paragraph{Verification of \eqref{21april16}.}
By rescaling it suffices to consider  $Q_0 = [0,1]\times [0,1]. $
For  $Q, Q' \in \cG $ with
$ |Q|= |Q'|$ and $Q\ne Q' $ we have $\la g_Q, g_{Q'} \ra= 0.$ Hence 
 the left hand side of \eqref{21april16} equals
\begin{equation} 
\label{21april18}
\sum_{Q \in \cG } \la g_Q, g_Q \ra
+ 
2\sum_{\{Q,Q' \in \cG :\, |Q| <| Q'|  \}} \la g_Q, g_{Q'} \ra.
\end{equation}
In view of \eqref{21april18} we aim at 
estimates for the entries of the Gram matrix 
$\la g_Q, g_{Q'}\ra .$ 

We first treat the diagonal terms of the
 Gram matrix.  
A  direct calculation gives 
$\la g_Q, g_{Q}\ra  =  \e |Q|/4,$ hence by \eqref{21april13}
\begin{equation} 
\label{21april21}
\sum_{Q \in \cG } \la g_Q, g_Q \ra \le  C  .
\end{equation}
Next we turn to estimating the off diagonal terms.
Consider  $Q,Q' \in \cG $ such  that  $|Q| <| Q'|. $
Write $Q = I \times J $ and $Q' = I' \times J'. $
Note, first if $\dist (Q,Q') \ge 2 \diam ( Q')$ then 
$\la g_Q, g_{Q'} \ra = 0.$
Hence it remains to consider the case   $\dist (Q,Q') \le 2 \diam ( Q') .$
Let $l_I$ denote the left endpoint of $I .$
The Lipschitz estimate $ \Lip(s_{I'}) \le C|I'|^{-1}$
and that $\int|d_{J}(x_2)| dx_2 \le \e |J| $
imply that 
\begin{equation} 
\label{21april19}
\begin{aligned}
|\la g_Q, g_{Q'} \ra |& = \left|\int ( s_{I'}(x_1) - s_{I'}(l_I))s_I(x_1)
                           dx_1\right|
                         \cdot \left|\int d_{J'}(x_2)d_{J}(x_2) dx_2 \right| \\
                       &\le C\frac{|I|}{|I'|}|I| \int|d_{J}(x_2)| dx_2 \\
                        &\le  \e C \frac{|I|}{|I'|} |Q|   .
\end{aligned}
\end{equation}
Since  $ Q = I \times J  \in \cG$ there exists $k \in \bN $ 
so that $ |I| = \e ^{2k} . $
Hence for $Q'= I'\times J'   \in \cG$ with $|Q'| >| Q| $ there exists
 $k' \in \bN $ with  $k'\le k -1 $ so that 
$ |I'| = \e ^{2k'} , $
and  $|I| / |I'| = \e ^{2k-2k' }.$
Note that for each $ Q \in \cG $ the cardinality of the set
$$  \{Q' \in \cG :\,\, |Q| <| Q'|,  \,\, \la g_Q, g_{Q'}\ra \ne 0 \}$$
is bounded by $C_1 ,$ say. 
Consequently in the double sum appearing on the left hand side of 
 \eqref{21april20}, for each $ Q $ only $C_1$ cubes  $Q'$ give a 
contribution.  
Thus by 
\eqref{21april19} 
\begin{equation} \label{21april20}
\begin{aligned}
\sum_{\{Q,Q' \in \cG :\, |Q| <| Q'|  \}} |\la g_Q, g_{Q'}\ra|
& \le C   \e ^{2k+1} \sum_{k' = 1 }^{k-1} \e ^{ - 2k'} \sum_{Q \in \cG   }|Q| \\
& \le  C\e^3 \sum_{Q \in \cG   }  |Q|.
\end{aligned}
\end{equation} 
By \eqref{21april13} the last line in \eqref{21april20}
is bounded by $ C  \e^2 . $
Combining \eqref{21april21} and \eqref{21april20}  gives \eqref{21april16}.

\paragraph{Verification of \eqref{21april17}.}
Fix $Q, Q_0 \in \cG$ so that $|Q_0| <| Q| $ and
$\dist( Q, Q_0 ) \le C\diam ( Q ). $
Then
\begin{equation} 
\label{21april22}
\|1_{Q_0}( g_Q - m_{ Q_0} ( g_Q) )\|_2 \le 
C\Lip (g_Q)\diam ( Q_0 ) |Q_0|^{1/2}.  
\end{equation}
Moreover if  $Q, Q_0 \in \cG$ so that $|Q_0| <| Q| $ and
$\dist( Q, Q_0 ) \ge C\diam ( Q ), $
then
\begin{equation} 
\label{21april23}
\|1_{Q_0}( g_Q - m_{ Q_0} ( g_Q)) \|_2 = 0 .
\end{equation}
Note that  
 $\Lip (g_Q) \le C(\e \diam ( Q ))^{-1}.$
Since  $ Q, Q_0 \in \cG ,$ with $|Q_0| <| Q| ,$ 
there exists $k,k_0 \in \bN ,$ with  $k\le k_0 -1 $
so that $\diam ( Q_0 ) = \sqrt{2} \cdot \e ^{2k_0}  $
and $\diam ( Q ) = \sqrt{2} \cdot \e ^{2k} . $
The cardinality of 
$$\{ Q\in \cG : \diam ( Q ) = \sqrt{2} \cdot \e ^{2k} ,\,\,\dist( Q, Q_0 ) \le C\sqrt{2} \cdot \e ^{2k} \} $$
is bounded by a constant $C . $ 
Hence by \eqref{21april22} and \eqref{21april23},
$$
\sum_{\{Q \in \cG ,\,|Q| >| Q_0 |\}}
\|1_{Q_0}( g_Q - m_{ Q_0} ( g_Q)) \|_2 
  \le C\e|Q_0|^{1/2}.
$$
Thus we verified \eqref{21april17}.
\endproof
We emphasize that the above upper bound on $\|f_\e\|_p $ 
works when the test functions $g = s\otimes d$ and its rescalings 
$g_Q = s_I \otimes d_J $ are defined with Lipschitz functions 
$A , B $ satisfying  \eqref{4sept075}, that is,  $\supp A \sbe [0,1] , \, \int A = 0 $ and 
$\supp B \sbe [-1,1] . $ We next impose furthermore that 
\begin{equation}
\label{4sept079}
{\rm A'\,\, is\,\, Lipschitz\,\, and } \int B = 0 . 
\end{equation}
 \begin{prop}\label{26april1}
Let $f_\e $ be defined by \eqref{21april12},
assume that \eqref{4sept079} and \eqref{4sept075} hold. 
Then for 
 $1 <p < \infty , $ 
$$\|R_1( f_\e)\|_{p} \le C_p \e   . $$
\end{prop}
\proof
The Fourier multipliers of the Riesz transforms
$R_1$ respectivley $R_2 $ are $\xi_1/|\xi|$ and 
 $\xi_2/|\xi| .$ Hence using \eqref{4sept079} 
for $g_Q = s_I\otimes d_J$ we have the identity
\begin{equation}\label{11okt071}
 R_1 ( g_Q) = R_2 (\pa _1 \bE_2 g_Q ) ,
\end{equation}
where $\pa _1 $ is differentiation with respect to the variable $x_1$ and 
$\bE_2 g_Q( x_1, x_2) = \int_{-\infty}^{x_2} g_Q(x_1, s)ds .$ 
Define now 
$$ \tilde s(x_1) = A'(x_1) ,
\quad \tilde d (x_2) =  C(x_2 /\e),
\quad  C(t) = \int_{-\infty}^{t} B(s) ds . $$
Let  $ \tilde s_I , \tilde d_J $ be  obtained from 
$\tilde s(x_1), \tilde d (x_2)$ by rescaling,
$$
\tilde s_I(x_1) = \tilde s(\frac{ x_1 - l_I}{|I|}), 
\quad
\tilde d_J(x_2) = \tilde d(\frac{ x_2 - l_J}{|J|}), 
$$ 
where  $l_I , l_J $ denote the left endpoint of 
$I$ respectively $J . $
Then with $\tilde g_Q = \tilde s_I\otimes \tilde d_J$
the identity \eqref{11okt071} assumes the following form,
\begin{equation}
\label{4sept0711}
R_1 (g_Q) = \e R_2(\tilde g_Q) .
\end{equation}

By \eqref{4sept079}
the Lipschitz functions
$ A' , C $ satisfy 
\eqref{4sept075}.
Hence Proposition~\ref{21april14} implies that 
$\tilde f_\e = \sum_{ Q \in  \cG } \tilde g_Q$
satisfies the $L^p$ estimate
$$ \| \tilde f_\e \|_p \le C_p . $$
By \eqref{4sept0711} we have 
$ R_1 (f_\e ) = \e R_2 ( \tilde f_\e ) .$ 
Hence the $L^p$ boundedness of the Riesz transforms yields 
$$ 
\begin{aligned} 
\| R_1 (f_\e )\|_p & \le \e \| R_2 (\tilde f_\e )\|_p\\
                   & \le C_p \e  \| \tilde f_\e \|_p\\
                   & \le   C_p \e.
\end{aligned}
$$
\endproof
We remark that the proof given above containd the following estimates
estimates that we will use again later.
For  $g = s \otimes d $ and $\tilde g = \tilde s \otimes \tilde d , $
\begin{equation}\label{15juli3}
 \begin{aligned}
\|R_1(g)\|_p    &= \e \| R_2 ( \tilde g)\|_p   \\
                & \le \e C_p \| \tilde g\|_p \\
                &  \le C_p\e^{1+1/p}. 
\end{aligned}
\end{equation}
\subsection{Lower bound for $\|P(f_\e )\|_p , \,p \ge 2 .$}
\label{lower}
We first  specialize once more the class of Lipschitz functions
$A,B$ we use to define  
$$
s(x_1) = A(x_1) , \quad d(x_2) = B(x_2 /\e )  $$
$$ g = s\otimes d \quad{\rm and} \quad   f_\e = \sum_{ Q \in  \cG } g_Q.
$$
We simply take now 
$$
B(x_2) =\begin{cases} 
           \sin(\pi x_2  ) & x_2 \in [-1,1];\\
           0             & x_2 \in \bR\sm  [-1,1].
        \end{cases}
$$
and choose $A$ to be smooth,  so that 
$\supp A \sbe [ 0,1] ,\,\int A = 0 $
and 
$$
\int_0^1 A(x_1) h_{[0,1]}(x_1) dx_1 =
\int_0^1 \sin( 2 \pi x_1 )  h_{[0,1]}(x_1) dx_1 .
$$
The following list of identities 
relates  the Haar functions $\{h_Q^{(1,0)}\}  $ to  
the test functions
 $\{g_Q\}.  $ 
\begin{enumerate}
\item
The scalar products $\la g_Q, h_Q^{(1,0)}\ra$ and $\la g_Q ,  g_Q \ra$ are as 
follows,
\begin{equation}
\label{21april52}
 \int g_Q(x) h_Q^{(1,0)}(x) dx = \e \frac {4 |Q|}{\pi^2} 
\quad\text{and}\quad
\int g_Q(x) g_Q(x) dx= \e \frac { |Q|}{4} . 
\end{equation}
\item
Let $Q' = I \times J' , $ be a dyadic square 
 where $J' $ is the dyadic interval 
adjacent to $J$ so that the right endpoint of $J$
is the left endpoint of $J'.$ Then
\begin{equation}
\label{21april50} 
\begin{aligned}
\int g_{Q'}(x) h_Q^{(1,0)}(x) dx &= -\int g_Q(x) h_Q^{(1,0)}(x) dx \\
&= -\e \frac {4 |Q|}{\pi^2}. 
\end{aligned}
\end{equation}
\item For all  choices of $Q' = I \times J'$ with $|J'|=|J|$
and $\dist(J,J') \ge |J| $ we have   
\begin{equation}
\label{21april51} 
\int g_{Q'}(x) h_Q^{(1,0)}(x) dx= 0.
\end{equation}
\item If $Q, Q' \in \cS $ so that $|Q'| < |Q| $ then 
\begin{equation}
\label{18sept06} 
 \int g_{Q'}(x) h_Q^{(1,0)} (x) = 0 . 
\end{equation}
\end{enumerate}

We consider  $ p \ge 2.$ Since  $P(f_\e )$ is compactly supported,
lower $L^p $ estimates  for  $P(f_\e )$ result  from lower $L^2 $ estimates. 
We obtain the latter by exploiting again 
the fact that $\{ g_Q : Q \in \cG \} $ 
is an ``almost orthogonal'' family of functions.

\begin{prop}\label{31april30}
Let $f_\e $ be defined by \eqref{21april12}.
The support of $P(f_\e) $ is contained in 
$ [-1,1]\times [-1,1] $ and
\begin{equation} 
\label{21april31} 
 \|P(f_\e)\|_{2} \ge c\e^{1/2} .
\end{equation}
Hence for $p\ge 2, $ $  \|P(f_\e)\|_{p} \ge c\e^{1/2} .$ 
\end{prop}
\proof
By Bessel's inequality, 
\begin{equation} 
\label{21april32} 
 \sum_{Q \in \cG  } \la f_\e , h_Q^{(1,0)}\ra ^2 
|Q|^{-1} \le  \|P(f_\e)\|_{2} ^2.
 \end{equation}
Using  \eqref{21april32} and \eqref{21april13}  
we prove below that  
\eqref{21april31} follows from the following lower estimate for the Haar 
coefficients
\begin{equation} 
\label{21april33}
|\la f_\e ,h_Q^{(1,0)}\ra| \ge c\e |Q| \quad\text{for}\quad Q \in \cG .
 \end{equation}
To prove \eqref{21april33},
fix   a dyadic square 
$Q = I \times J $ with $ Q \in \cG . $
Write  the Haar coefficient as 
\begin{equation} 
\label{21april330}
\la f_\e ,h_Q^{(1,0)}\ra = \la  g_Q , h_Q^{(1,0)}\ra + 
\sum_{Q' \in \cG \sm \{ Q \}  } \la  g_{Q'}, h_Q^{(1,0)} \ra .
\end{equation}
Recall \eqref{21april52}
asserting that 
$$\la  g_Q , h_Q^{(1,0)}\ra = \e 4|Q|/\pi^2 . $$
Next we show that the off diagonal terms in \eqref{21april330}
are negligible compared to $\la  g_Q , h_Q^{(1,0)}\ra.$
We claim,
\begin{equation} 
\label{21april34}
\sum_{Q' \in \cG \sm \{ Q \}  } |\la  g_{Q'}, h_Q^{(1,0)} \ra| 
\le C \e^2 |Q| . 
\end{equation}
The first step in the verification of the claim consists in observing that 
the only contribution to 
 \eqref{21april34} comes from the index set $\{Q' \in \cG\sm \{ Q \}:\,|Q'| >| Q |  \}. $ Indeed, if    $Q' \in \cG,$   $Q' \ne Q$ and   $|Q'| = | Q |$
then  \eqref{21april51} in combination with 
\eqref{21april1000} implies  that  $\la  g_{Q'}, h_Q^{(1,0)} \ra =0. $
Also by \eqref{18sept06}
for   $Q' \in \cG$ and 
$|Q'| < | Q |  $ we 
have $\la  g_{Q'}, h_Q^{(1,0)} \ra =0. $

Next we provide an estimate for  the contribution to \eqref{21april34}
coming from $\{Q' \in \cG\sm \{ Q \}:\,|Q'| >| Q |  \}. $
Choose $ k \in \bN $ so that $| Q | = \e^{4k}  $
and let $ k' \in \bN $ satisfy 
$ k'<k .$ There exists at most one square
$Q' \in   \cG $ satisfying 
$$ |Q'| = \e^{4k'}
\quad\text{and}\quad 
\la  g_{Q'}, h_Q^{(1,0)} \ra \ne 0 . $$
Next fix  
$Q' = I' \times J' $ with $ |Q'| = \e^{4k'}$
and  $ k'<k .$ 
Write $Q = I \times J $ and $Q' = I' \times J'. $
Let $l_I$ denote the left endpoint of $I .$
Recall that $\Lip( s_{I'}) \le C|I'|^{-1}$
and  $\int|d_{J}(x_2)| dx_2 \le C |J|. $
Hence,
\begin{equation} 
\label{21april35}
\begin{aligned}
|\la g_{Q'},  h_Q^{(1,0)} \ra |
& = \left|\int ( s_{I'}(x_1) - s_{I'}(l_I))h_I(x_1)
                           dx_1\right|
                         \cdot \left|\int_J d_{J'}(x_2) dx_2 \right| \\
                        &\le C  |I|\cdot|I'|^{-1} |Q|\\
                        & = C\e^{2k - 2k'}  |Q|   .
\end{aligned}
\end{equation}
By definition of  $g_{Q'}$ and $ h_Q^{(1,0)}$ if 
 $|Q'| >| Q |$ and $\la  g_{Q'}, h_Q^{(1,0)} \ra \ne 0 $ then
$ \dist( Q', Q ) \le C \diam (Q' ) . $  
It now follows  from \eqref{21april35} that for any $Q \in \cG , $
\begin{equation} 
\label{21april36}
\begin{aligned}
\sum_{\{Q' \in \cG ,\, |Q'| >| Q |  \}}
 |\la  g_{Q'}, h_Q^{(1,0)} \ra | 
& \le C |Q|  \e^{  2k} \sum_{k' = 1 } ^{k-1} \e^{-2k'}   \\
& \le C \e^2 |Q| .
\end{aligned}
\end{equation}
Thus by \eqref{21april36} we verified the claim  \eqref{21april34}.
Hence  we have \eqref{21april33}.
It remains to show how  the coefficient estimates 
\eqref{21april33} imply the norm inequality of \eqref{21april31}.
Using first \eqref{21april32}   then \eqref{21april33} and
\eqref{21april13} we obtain 
$$
\begin{aligned}
\|P(f_\e) \|_2^2 
                 &\ge c \e^2 \sum_{Q\in \cG } |Q| \\
                 & \ge c\e . 
\end{aligned}
$$

\endproof

\subsection{The proof of theorem~\ref{th200} .
}
\label{proof}
We choose 
Lipschitz functions $A, B $ 
with specification of the 
previous sub-section and define  
testing functions $ g = s_{[0,1]}\otimes d_{[0,1]} ,$
 $f_\e $
as above.

Consider first  the estimate \eqref{14april1111} of Theorem~\ref{th200}.
Let $ 1 <  p \le 2 . $  Fix $  \eta > 0 .$ Let $ g = s_{[0,1]}\otimes d_{[0,1]} $
be defined by \eqref{15juli1}. 
Since $g$ is bounded and supported in $[0,1]\times [-\e,\e] , $ 
we have 
\begin{equation}\label{15juli2}
 \|g\|_p \le C\e^{1/p} .
\end{equation}

Next observe that for  the square function $\bS (P(g))$
we have the obvious  estimate 
$\bS (P(g)) \ge  | \la g , h^{(1,0)}_{[0,1[\times[0,1[}\ra| .$
Next recall that 
$\|P(g)\|_p \sim  \|\bS (P(g) ) \|_p$ 
hence $\|P(g)\|_p \ge c | \la g , h^{(1,0)}_{[0,1[\times[0,1[}\ra| .$
By \eqref{21april52},  we have  
$\la g , h^{(1,0)}_{[0,1[\times[0,1[}\ra = 4\e / \pi^2 , $
hence
\begin{equation}\label{15juli4}
 ||P(g)||_p \ge c\e .
\end{equation}
By \eqref{15juli2} and \eqref{15juli3}
\begin{equation}\label{15juli5}
 \|g\|_p^{1/p-\eta}\|R_1(g) \|_p^{1/q+\eta}\le C \e^{1+\eta} .
\end{equation}
Combining \eqref{15juli4} and \eqref{15juli5} yields
$$\frac{||P(g)||_p }{ \|g\|_p^{1/p-\eta}\|R_1(g)\|_p^{1/q+\eta}}
\ge c\e^{-\eta} .
  $$
 Since $\eta > 0 $ is fixed and $\e >0 $ is arbitrarily small 
we verified \eqref{14april1111}.

Next we turn to the case 
$p\ge 2.$
The test function  $f_\e $ is defined by  \eqref{21april12}.
Proposition~\ref{21april14} and Proposition~\ref{26april1}
 give the upper bounds
$$ \|f_\e\|_p \le C_p \quad\text{and}\quad\|R_1(f_\e)\|_p \le 
C_p\e . $$ 
Hence for $ \eta >0 $
$$ \|f_\e\|_p^{1/2 -\eta} \|R_1(f_\e)\|_p^{1/2 +\eta}
\le C_p\e^{1/2 + \eta} . $$
By Proposition~\ref{31april30} we have the lower estimate 
$$\|P(f_\e)\|_p \ge c_p \e^{1/2}.$$
so that 
$$ \frac{\|P(f_\e)\|_p}{ \|f_\e\|_p^{1/2 -\eta} \|R_1(f_\e)\|_p^{1/2 +\eta}}
\ge c_p \e^{- \eta} . $$
\endproof
\paragraph{Acknowledgement.} 
During the preparation of this paper J. Lee and P. M\"uller 
visited repeatedly the Max Planck Insitute for Mathematics in the Sciences in 
Leipzig. It is our pleasure to thank this institution for its  
hospitality and for providing us with excellent working conditions.
\nocite{jpw1}
 \bibliographystyle{abbrv}
\bibliography{compensated}

\noindent
{\bf  AMS Subject classification: 49J45 , 42C15, 35B35}\\

\noindent
{\bf  Addresses:}\\

\noindent
Jihoon Lee\\
Department of Mathematics\\
Sungkyunkwan University\\ 
Suwon, Korea\\
E-mail:jihoonlee@skku.edu
\vskip 0.3cm
\noindent
 Paul F.X. M\"uller\\
Institut f\"ur Analysis \\
J. Kepler Universit\"at\\
 A-4040 Linz \\
E-mail: pfxm@bayou.uni-linz.ac.at
\vskip 0.3cm
\noindent 
Stefan M\"uller\\
Max Planck Institute for Mathematics in the Sciences\\
Inselstr. 22--26\\
D 04103 Leipzig Germany\\
E-mail: sm@mis.mpg.de
\vskip 0.5cm
\noindent
\end{document}